\documentclass{colt2016} 

\title[Optimal Best Arm Identification with Fixed Confidence]{Optimal Best Arm Identification with Fixed Confidence}

\usepackage{times}
\usepackage{macrosBonus}

 \coltauthor{\Name{Aur\'elien Garivier} \Email{aurelien.garivier@math.univ-toulouse.fr}\\
 \addr Institut de Math\'ematiques de Toulouse; UMR5219\\
 Universit\'e de Toulouse; CNRS\\
 UPS IMT, F-31062 Toulouse Cedex 9, France
 \AND
 \Name{Emilie Kaufmann} \Email{emilie.kaufmann@inria.fr}\\
 \addr Univ. Lille, CNRS, UMR 9189 \\
CRIStAL - Centre de Recherche en Informatique Signal et Automatique de Lille\\
 F-59000 Lille, France
 }

\begin{document}

\maketitle

\begin{abstract}
We give a complete characterization of the complexity of best-arm identification in one-parameter bandit problems. We prove a new, tight lower bound on the sample complexity. We propose the `Track-and-Stop' strategy, which we prove to be asymptotically optimal. It consists in a new sampling rule (which tracks the optimal proportions of arm draws highlighted by the lower bound) and in a stopping rule named after Chernoff, for which we give a new analysis. 
\end{abstract}


\section{Introduction}

A multi-armed bandit model is a paradigmatic framework of sequential statistics made of $K$ probability distributions $\nu_1,\dots,\nu_K$ with respective means $\mu_1,\dots,\mu_K$: at every time step $t=1,2,\dots$ one \emph{arm} $A_t\in\Arms=\{1,\dots,K\}$ is chosen and a new, independent \emph{reward} $X_t$ is drawn from $\nu_{A_t}$.
Introduced in the 1930s with motivations originally from clinical trials, bandit models have raised a large interest recently as relevant models for interactive learning schemes or recommender systems.
A large part of these works consisted in defining efficient strategies for maximizing the expected cumulated reward $\bE[X_1+\dots+X_t]$; see for instance~\cite{Bubeck:Survey12} for a survey. 
A good understanding of this simple model has allowed for efficient strategies in much more elaborate settings, for example including side information~\citep{YadLinear11,AGContext13}, infinitely many arms~\citep{SrinivasGPUCB,Bubeck11Xarmed}, or for the search of optimal strategies in games~\citep{SurveyRemiMCTS}, to name just a few. 

In some of these applications, the real objective is not to maximize the cumulated reward, but rather to identify the arm that yields the largest mean reward $\mu^* = \max_{1\leq a\leq K}\mu_a$, as fast and accurately as possible, regardless of the number of bad arm draws. Let ${\cal F}_t=\sigma(X_1,\dots,X_t)$ be the sigma-field generated by the observations up to time $t$.  A strategy is then defined by:
\begin{itemize}[noitemsep,topsep=0pt,parsep=0pt,partopsep=0pt]
\item a \emph{sampling rule} $(A_t)_t$, where $A_t$ is ${\cal F}_{t-1}$-measurable;
\item a \emph{stopping rule} $\tau$, which is a stopping time with respect to ${\cal F}_t$;
\item and a ${\cal F}_{\tau}$-measurable \emph{decision rule} $\hat{a}_\tau$.
\end{itemize}
The goal is to guarantee that $\hat{a}_\tau\in\text{argmax}\, \mu_a$ with the highest possible probability while minimizing the number $\tau$ of draws.
Two settings have been considered in the literature. In the \emph{fixed-budget setting}, 
the number of draws is fixed in advance, and one aims at minimizing the probability of error $\bP(\hat{a}_\tau\notin\text{argmax}\, \mu_a)$. In the \emph{fixed-confidence setting} 
a maximal risk $\delta$ is fixed, and one looks for a strategy guaranteeing that $\bP(\hat{a}_\tau\notin\text{argmax}\, \mu_a) \leq \delta$ (such a strategy is called \emph{$\delta$-PAC}) while minimizing the \emph{sample complexity} $\bE[\tau]$.

The aim of this paper is to propose an analysis of the best arm identification problem in the fixed-confidence setting. For the sake of clarity and simplicity, we suppose that there is a single arm with highest expectation, and without loss of generality that $\mu_1>\mu_2\geq\dots\geq \mu_K$.
We also focus on the simple case where the distributions are parameterized by their means, as in one-parameter exponential families, and we  index probabilities and expectations by the parameter ${\bm \mu}=(\mu_1,\dots,\mu_K)$. The Kullback-Leibler divergence of two distributions of means $\mu_1$ and $\mu_2$ is a function $d : (\mu_1,\mu_2)\to\R^+$. Two cases are of particular interest: the Gaussian case, with $d(x,y)=(x-y)^2/(2\sigma^2)$, and binary rewards with $d(x,y)=\kl(x,y):=x\log(x/y)+(1-x)\log((1-x)/(1-y))$.

Several strategies have been proposed to minimize $\bE_{\bm \mu}[\tau]$. While racing strategies are based on successive eliminations of apparently sub-optimal arms (\cite{EvenDaral06,COLT13}), another family of strategies exploits the use of upper and lower confidence bounds on the means of the arms (\cite{Shivaramal12,Gabillon:al12,COLT13,Jamiesonal14LILUCB}). 
They reflect some aspects of the difficulty of the problem, but are not proved to satisfy any optimality property\footnote{Optimality is mentioned in several articles, with different and sometimes weak meanings (minimax, rate-optimal,...). In our view, BAI algorithms for which there exists a model with a sample complexity bounded, up to a multiplicative constant, by some quantity related to some lower bound, may not be called  optimal.}. In particular, there was still a gap between the lower bounds, involving complexity terms reflecting only partially the structure of the problem, and the upper bounds on $\bE_{\bm \mu}[\tau]$ for these particular algorithms, even from an asymptotic point of view. For the particular case $K=2$, this gap was closed in \cite{COLT14}. The tools used there, however, where specific to the two-arm case and cannot be extended easily.

The first result of this paper is a tight, non-asymptotic lower bound on $\bE_{\bm \mu}[\tau]$. 
This bound involves a `characteristic time' for the problem, depending on the parameters of the arms, which does not take a simple form like for example a sum of squared inverse gaps. Instead, it appears as the solution of an optimization problem, in the spirit of the bounds given by~\cite{GravesLai97} in the context of regret minimization.
We give a brief analysis of this optimization problem, and we provide an efficient numerical solution as well as elements of interpretation. 

The second contribution is a new $\delta$-PAC algorithm that asymptotically achieves this lower bound, that we call the Track-and-Stop strategy.
In a nutshell, the idea is to sample so as to \emph{equalize the probability of all possible wrong decisions}, and to stop as soon as possible. The stopping rule, which we name after Chernoff, can be interpreted in three equivalent ways: in statistical terms, as a generalized likelihood ratio test; in information-theoretic terms, as an application of the Minimal Description Length principle; and in terms of optimal transportation, in light of the lower bound. The sampling rule is a by-product of the lower bound analysis, which reveals the existence of optimal proportions of draws for each arm. By estimating and tracking these proportions, our algorithm asymptotically reaches the optimal sample complexity as $\delta$ goes to~$0$.

The paper is organized as follows. 
In Section~\ref{sec:lb}, the lower bound is given with a short proof.
Section~\ref{sec:algo} contains a commented description of our stopping and sampling rules.
The analysis of the resulting algorithm is sketched in Sections~\ref{sec:PACanalysis} (validity of the stopping rule) and~\ref{sec:SCAnalysis} (sample complexity analysis), establishing its asymptotic optimality.
Section~\ref{sec:experiments} contains practical comments and results of numerical experiments, in order to show the  efficiency of our strategy even for moderate values of $\delta$. We also briefly comment on the gain over racing strategies, which can be explained in light of Theorem~\ref{cor:GeneBAI}. Most proofs and technical details are postponed to the Appendix.

\section{Lower Bounds on the sample complexity}\label{sec:lb}
The pioneering work of ~\cite{LaiRobbins85bandits} has popularized the use of changes of distributions to show problem-dependent lower bounds in bandit problems: the idea is to move the parameters of the arms until a completely different behavior of the algorithm is expected on this alternative bandit model. The cost of such a transportation is induced by the deviations of the arm distributions: by choosing the most economical move, one can prove that the alternative behavior is not too rare in the original model. Recently, \cite{COLT14} and \cite{Combes14Unimodal} have independently introduced a new way of writing such a change of measure, which relies on a transportation lemma that encapsulates the change of measure and permits to use it at a higher level. 
Here we go one step further by combining \emph{several changes of measures at the same time}, in the spirit of~\cite{GravesLai97}. This allows us to prove a non-asymptotic lower bound on the sample complexity valid for any $\delta$-PAC algorithm on any bandit model with a unique optimal arm. We present this result in the particular case of arms that belong to a canonical exponential family,
\[\cP = \left\{(\nu_\theta)_{\theta \in \Theta} : \frac{d\nu_\theta}{d\xi} = \exp\big(\theta x - b(\theta)\big)  \right\},\]
where ${\Theta}\subset{\R}$, $\xi$ is some reference measure on $\R$ and $b : \Theta \rightarrow \R$ is a convex, twice differentiable function. A distribution $\nu_\theta\in\cP$ can be parameterized by its mean $\dot{b}(\theta)$, and for every $\mu \in \dot{b}(\Theta)$ we denote by $\nu^\mu$ be the unique distribution in $\cP$ with expectation $\mu$. Unfamiliar readers may simply think of Bernoulli laws, or Gaussian distributions with known variance. As explained in~\citet[see also references therein]{KLUCBJournal}, the Kullback-Leibler divergence from $\nu_\theta$ to $\nu_{\theta'}$ induces a divergence function $d$ on $\dot{b}(\Theta)$ defined, if $\dot{b}(\theta)=\mu$ and $\dot{b}(\theta')=\mu'$,~by
\[
d(\mu,\mu') = \text{KL}(\nu^{\mu},\nu^{\mu'}) = \text{KL}(\nu_\theta,\nu_{\theta'})=b(\theta')-b(\theta)-\dot{b}(\theta)(\theta'-\theta) \;.\]
With some abuse of notation, an exponential family bandit model $\nu=(\nu^{\mu_1},\dots,\nu^{\mu_K})$ is identified with the means of its arms  $\bm \mu = (\mu_1,\dots,\mu_K)$.

\subsection{General Lower Bound}
Denote by $\cS$ a set of exponential bandit models such that each bandit model $\bm{\mu} =(\mu_1,\dots,\mu_K)$ in $\cS$ has a unique optimal arm: for each $\bm{\mu}\in \cS$, there exists an arm $a^*(\bm \mu)$ such that $\mu_{a^*(\bm \mu)} > \max \{ \mu_a : a \neq a^*(\bm \mu)\}$.
Fixed-confidence strategies depend on the risk level and we subscript stopping rules by $\delta$. A strategy is called $\delta$-PAC if for every $\bm \mu \in \cS$, $\bP_{\bm \mu}(\tau_\delta < \infty) = 1$ and  $\bP_{\bm \mu}(\hat{a}_{\tau_\delta} \neq a^*) \leq \delta$.
We introduce \[\Alt(\bm \mu) := \{ \bm\lambda \in \cS : a^*(\bm\lambda) \neq a^*(\bm \mu)\}\;,\]
the set of problems where the optimal arm is not the same as in $\bm\mu$, and $\Sigma_K=\{\omega\in\R^k_+ : \omega_1+\dots+\omega_K=1\}$ the set of probability distributions on $\Arms$.
\begin{theorem}\label{cor:GeneBAI} Let $\delta \in (0,1)$. 
For any $\delta$-PAC strategy and any bandit model $\bm\mu\in\cS$,  
\[\bE_{\bm \mu}[\tau_\delta] \geq T^*(\bm\mu)\,\kl(\delta,1-\delta), \]
where 
\begin{equation}
T^*(\bm \mu)^{-1} := \sup_{w \in \Sigma_K} \inf_{\bm \lambda \in \Alt(\bm\mu)} \left( \sum_{a=1}^K w_a d(\mu_a,\lambda_a)\right)\;.\label{equ:Original} 
\end{equation}
\end{theorem}

\begin{remark} As $\mathrm{kl}(\delta,1-\delta) \sim \log(1/\delta)$ when $\delta$ goes to zero, Theorem~\ref{cor:GeneBAI} yields the asymptotic lower bound
\[\liminf_{\delta \rightarrow 0} \frac{\bE_{\bm \mu}[\tau_\delta]}{\log(1/\delta)} \geq T^*(\bm \mu).\]
A non-asymptotic version can be obtained for example from the inequality $\kl(\delta,1-\delta) \geq \log(1/(2.4\delta))$ that holds for all $\delta\in(0,1)$, given in~\cite{JMLR15}.
\end{remark}
We will see that the supremum in Equation~\eqref{equ:Original} is indeed a maximum, and we call   
\[w^*(\bm\mu): = \argmax{w \in \Sigma_K} \ \inf_{\bm \lambda \in \Alt(\bm\mu)} \left( \sum_{a=1}^K w_a d(\mu_a,\lambda_a)\right)\]
the corresponding distribution on the arms. 
The proof of Theorem~\ref{cor:GeneBAI} shows that $w^*$ is the  proportion of arm draws of any strategy matching this lower bound.

The particular case where $\cS$ is the class of bandit models with Poisson rewards in which all suboptimal arms are equal is considered in the very insightful paper by~\cite{Rajesh15Oddball}, where a closed-form formula is given for both $T^*(\bm \mu)$ and $w^*(\bm\mu)$. In this paper, we consider best arm identification in all possible bandit models with a single best arm, and in the rest of the paper we fix 
\begin{equation*}\cS := \left\{ \bm \mu : \exists a \in \Arms : \mu_{a} > \max_{i\neq a} \mu_i \right\}.\end{equation*}

\paragraph{Proof of Theorem~\ref{cor:GeneBAI}.}
Let $\delta \in (0,1)$, $\bm \mu \in \cS$, and consider a $\delta$-PAC strategy. 
For every $t\geq 1$, denote by $N_a(t)$ the (random) number of draws of arm $a$ up to time $t$.
The `transportation' Lemma 1 of~\cite{JMLR15} relates the expected number of draws of each arm and the Kullback-Leibler divergence of two bandit models with different optimal arms to the probability of error $\delta$: 
\begin{equation}
\forall \lambda \in \cS : a^*(\bm\lambda) \neq a^*(\bm\mu), \ \ \ \sum_{a=1}^K d(\mu_a,\lambda_a) \bE_{\bm \mu}[N_a(\tau_\delta)] \geq  \kl(\delta,1-\delta).
 \label{eq:BAI}
\end{equation}
Instead of choosing for each arm $a$ a specific instance of $\bm\lambda$ that yields a lower bound on $\bE_{\bm \mu}[N_a(\tau_\delta)]$, we combine here the inequalities given by all alternatives $\bm\lambda$: 
\begin{align*}
\kl(\delta,1-\delta) &\leq  \inf_{\lambda \in \Alt(\bm \mu)} \!\bE_{\bm\mu}[\tau_\delta] \left(\sum_{a=1}^{K} \frac{\bE_{\bm\mu}[N_a]}{\bE_{\bm\mu}[\tau_\delta]}d(\mu_a,\lambda_a)\right) =  \bE_{\bm\mu}[\tau_\delta]\inf_{\lambda \in \Alt(\bm \mu)} \left(\sum_{a=1}^{K} \frac{\bE_{\bm\mu}[N_a]}{\bE_{\bm\mu}[\tau_\delta]}d(\mu_a,\lambda_a)\right)\\
&\leq  \bE_{\bm\mu}[\tau_\delta] \sup_{w \in \Sigma_K}\!\inf_{\lambda \in \Alt(\bm \mu)}\! \left(\sum_{a=1}^{K} w_ad(\mu_a,\lambda_a)\right),
\end{align*}
as $\bE_{\bm\mu}[\tau_\delta] = \sum_{a=1}^K \bE_{\bm\mu}[N_a(\tau_\delta)]$.\qed 
In the last inequality, the strategy-dependent proportions of arm draws are replaced by their supremum so as to obtain a bound valid for any $\delta$-PAC algorithm; one can see that a strategy may reach this bound only if it meets the $w^*(\bm\mu)$.
To make this bound useful, it remains to study $T^*$ and $w^*$.

\subsection{About the Characteristic Time and the Optimal Proportions}
We study here the optimization problem~\eqref{equ:Original}, so as to better understand the function $T^*$ and $w^*$ (Proposition~\ref{prop:PropNu}), and in order to provide an efficient algorithm for computing first $w^*(\bm\mu)$ (Theorem~\ref{thm:ExplicitForm}), then $T^*(\bm\mu)$ (Lemma~\ref{lem:CDinterpretation}). The main ideas are outlined here, while all technical details are postponed to Appendix~\ref{proofs:Nu}. 
Simplifying $T^*(\bm \mu)$ requires the introduction of the following parameterized version of the Jensen-Shannon divergence (which corresponds to $\alpha=1/2$): for every $\alpha\in[0,1]$, let
\begin{equation}I_\alpha(\mu_1,\mu_2):=\alpha d\big(\mu_1,\alpha\mu_1 + (1-\alpha)\mu_2\big) + (1-\alpha)d\big(\mu_2,\alpha\mu_1 + (1-\alpha)\mu_2\big)\;.\label{def:Divergences}\end{equation}
The first step is, for any $w$, to identify the minimizer of the transportation cost: 
\begin{lemma}\label{lem:CDinterpretation} For every $w\in \Sigma_K$,
\[ \inf_{\bm\lambda \in \Alt(\bm \mu)} \left(\sum_{a=1}^{K} w_ad(\mu_a,\lambda_a)\right) = \min_{a \neq 1} \ (w_1 + w_a)I_{\frac{w_1}{w_1+w_a}}(\mu_1,\mu_a)\;.\]
\end{lemma}
It follows that 
\begin{eqnarray*}T^*(\bm\mu)^{-1} 
 & = & \sup_{w \in \Sigma_K} \min_{a \neq 1} \ (w_1 + w_a)I_{\frac{w_1}{w_1+w_a}}(\mu_1,\mu_a)\; ,\hbox{and} \\
w^*(\bm\mu) & = & \argmax{w \in \Sigma_K} \ \min_{a \neq 1} \ (w_1 + w_a)I_{\frac{w_1}{w_1+w_a}}(\mu_{1},\mu_a)\;.
\end{eqnarray*}
It is easy to see that, at the optimum, the quantities $(w_1 + w_a)I_{{w_1}/({w_1+w_a})}(\mu_1,\mu_a)$ are all equal.
\begin{lemma}\label{lem:allEqual}
For all $a,b\in\{2,\dots,K\}$, 
\[ (w^*_1 + w^*_a)I_{\frac{w^*_1}{w^*_1+w^*_a}}(\mu_1,\mu_a) =  (w^*_1 + w^*_b)I_{\frac{w^*_1}{w^*_1+w^*_b}}(\mu_1,\mu_b)\;. \]
\end{lemma}
This permits to obtain a more explicit formula for $w^*(\bm \mu)$ involving only a single real parameter. Indeed, for every $a \in\{ 2, \dots K\}$ let 
 \begin{equation}g_a(x)  =   (1+x) I_{\frac{1}{1+x}}(\mu_1,\mu_a) \;.
\label{def:InverseFunction}\end{equation}
The function $g_a$ is a strictly increasing one-to-one mapping from $[0,+\infty[$ onto $[0,d(\mu_1,\mu_a)[$. We define $x_a: [0,d(\mu_1,\mu_a)[ \to [0, + \infty[$ as its inverse function: $x_a(y)=g_a^{-1}(y)$.
Denoting $x_1$ the function constantly equal to $1$, one obtains the following characterization of $w^*(\bm\mu)$:
\begin{theorem}\label{thm:ExplicitForm} For every $a\in\Arms$, 
\begin{equation}\label{eq:nuofx}
w^*_a(\bm \mu) = \frac{x_a(y^*)}{\sum_{a=1}^Kx_a(y^*)}\;,\end{equation}
where $y^*$ is the unique solution of the equation $F_{\bm\mu}(y)=1$, and where
\begin{equation}\label{eq:relation} F_{\bm\mu}:y\mapsto \sum_{a=2}^K \frac{d\left(\mu_1,\frac{\mu_1 + x_a(y)\mu_a}{1+x_a(y)}\right)}{d\left(\mu_a,\frac{\mu_1 + x_a(y)\mu_a}{1+x_a(y)}\right)} \end{equation}
is a continuous, increasing function on $[0,d(\mu_1,\mu_2)[$ such that $F_{\bm\mu}(0)=0$ and $F_{\bm\mu}(y)\to \infty$ when $y\to d(\mu_1,\mu_2))$.
\end{theorem}
Thus, $w^*$ can be simply computed by applying (for example) the bisection method to a function whose evaluations requires the resolution of $K$ smooth scalar equations. By using efficient numerical solvers, we obtain a fast algorithm of complexity, roughly speaking, proportional to the number of arms. 
This characterization of $w^*(\bm\mu)$ also permits to obtain a few sanity-check properties, like for example:
\begin{proposition}\label{prop:PropNu}
 \begin{enumerate}
  \item For all $\bm\mu \in \cS$, for all $a$, $w_a^*(\bm\mu) \neq 0$. 
  \item $w^*$ is continuous in every $\bm \mu \in \cS$. 
  \item If $\mu_1>\mu_2 \geq \dots \geq \mu_K$, one has $w_2^*(\bm \mu) \geq \dots \geq w_K^*(\bm \mu)$.
 \end{enumerate}
\end{proposition}
Observe that one may have\footnote{This can happen when the right-deviations of $\nu^{\mu_2}$ are smaller than the left-deviations of $\nu^{\mu_1}$; for example, with Bernoulli arms of parameters $\bm\mu=(0.5,0.1, 0.02)$, $\nu^*(\bm\mu) \approx(39\%,42\%, 19\%)$.} $w_2>w_1$. In general, it is not possible to give closed-form formulas for $T^*(\bm \mu)$ and $w^*(\bm\mu)$. 
In particular, $T^*(\bm \mu)$ cannot be written as a sum over the arms of individual complexity terms as in previous works (\cite{MannorTsi04,JMLR15}). 
But the following particular cases can be mentioned.

\paragraph{Two-armed bandits.} 
For a two-armed bandit model $\bm \mu=(\mu_1,\mu_2)$, $T^*(\bm\mu)$ and $w^*(\bm \mu)$ can be computed algebraically. 
Lemma~\ref{lem:CDinterpretation} and the fact that $w_2=1-w_1$ imply that 
\begin{eqnarray*}T^*(\bm\mu)^{-1} 
 & = & \sup_{\alpha \in (0,1)} I_{\alpha}(\mu_1,\mu_a),
\end{eqnarray*}
Some algebra shows that the maximum is reached at $\alpha=\alpha_*(\mu_1,\mu_2)$ defined by the equation  $d(\mu_1,\mu_*) = d(\mu_2,\mu_*)$, where $\mu_*=\alpha_*(\mu_1,\mu_2)\mu_1 + (1-\alpha_*(\mu_1,\mu_2))\mu_2$.
The value of the maximum is then the `reversed' Chernoff information
$d_*(\mu_1,\mu_2):=d(\mu_1,\mu_*)=d(\mu_2,\mu_*)$.
This permits to recover the bound already given in~\cite{COLT14}:
\[\bE_{\bm \mu}[\tau_\delta] \geq \frac{\kl(\delta, 1-\delta)}{d_*(\mu_1,\mu_2)} \;.\]

\paragraph{The Gaussian case.}
When $d(x,y) = (x-y)^2/(2\sigma^2)$, $T^*(\mu)$ and $w^*(\bm \mu)$ can be computed by solving a rational equation. Indeed, Equation~\eqref{eq:relation} and Lemma~\ref{lem:allEqual} imply that 
\[\sum_{a=2}^K x_a^2 = 1 \ \ \ \ \text{and} \ \ \ \ \frac{x_a}{1+x_a} = \frac{\lambda}{(\mu_1-\mu_a)^2}, \ \ \text{thus} \ \ \sum_{a=2}^K\frac{\lambda^2}{\big((\mu_1-\mu_a)^2 - \lambda\big)^2} = 1\]
for some $\lambda \in (0, (\mu_1-\mu_2)^2)$. For $K=3$, $\lambda$ is the solution of a polynomial equation of degree 4 
and has therefore an (obscure) algebraic expression. 
The following inequalities, established in Appendix~\ref{proof:BoundGaussian} give a better insight of the order of magnitude of $T^*(\bm \mu)$: if $\Delta_1 = \Delta_2$ and $\Delta_a = \mu_1 - \mu_a$ for $a\geq 2$, then
\[\sum_{a=1}^K \frac{2\sigma^2}{\Delta_a^2} \leq T^*(\bm \mu) \leq 2 \sum_{a=1}^K \frac{2\sigma^2}{\Delta_a^2}.\]

\section{The Track-and-Stop Strategy}\label{sec:algo}
We now describe a new strategy which is the first (as far as we know) to asymptotically match the lower bound of Theorem~\ref{cor:GeneBAI}.
Denote by $\hat{\mu}(t) = (\hat{\mu}_1(t), \dots, \hat{\mu}_K(t))$ the current maximum likelihood estimate of $\bm\mu$ at time $t$: $\hat{\mu}_a(t)=N_a(t)^{-1}\sum_{s\leq t}X_s \ind\{A_s=a\}$.
As seen in Section~\ref{sec:lb}, a good sampling rule should respect the optimal proportions of arm draws given by $w^*(\bm \mu)$. There are several ways to ensure this, and we present two of them in Section~\ref{subsec:Sampling}.
A good stopping rule should determine the earliest moment when sufficient statistical evidence has been gathered for identifying the best arm: we propose one (with several interpretations) in Section~\ref{subsec:Stopping}, showing in Section~\ref{sec:PACanalysis} how to tune it in order to ensure the $\delta$-PAC property.
As for the decision rule, we simply choose $\hat{a}_{\tau_\delta} = \argmax{a\in\Arms} \ \hat{\mu}_a(\tau_\delta)$.
The optimality of the Track-and-Stop strategy is shown in Section~\ref{sec:SCAnalysis}.
\subsection{Sampling Rule: Tracking the Optimal Proportions}
\label{subsec:Sampling}
The first idea for matching the proportions $w^*(\bm \mu)$ is to track the plug-in estimates $w^*(\hat{\mu}(t))$. In bandit settings, using plug-in estimates is always hazardous, because bad initial estimate may lead to abandon an arm and to prevent further observations that would correct the initial error. Indeed, one may see (both theoretically and in numerical experiments) that a naive plug-in sampling rule sometimes fails. But there is a very simple fix, which consists in forcing sufficient exploration of each arm to ensure a (sufficiently fast) convergence of $\hat{\mu}(t)$. 

The shortest way to do this (which we term \emph{C-Tracking}) is to slightly alter the optimization solution: for every $\epsilon\in(0,1/K]$, let $w^\epsilon(\bm \mu)$ be a $L^\infty$ projection of $w^*(\bm \mu)$ onto $ \Sigma^\epsilon_K=\big\{(w_1,\dots,w_K)\in [\epsilon, 1] : w_1+\dots+w_K=1\big\}$.
Choosing $\epsilon_t = (K^2+t)^{-1/2}/2$ and 

\vspace{-0.4cm}\[A_{t+1} \in \argmax{1\leq a\leq K} \sum_{s=0}^{t} w^{\epsilon_s}_a(\hat{\bm\mu}(s)) - {N_a}(t)\;,\]

\vspace{-0.3cm}\noindent we prove in Appendix~\ref{proof:Tracking} that: 
\begin{lemma}\label{lem:TrackingC}
For all $t\geq 1$ and $a\in\Arms$, the C-Tracking rule ensures that $N_a(t)\geq\sqrt{t+K^2}-2K$ and that
\[\max_{1\leq a\leq K} \left|N_a(t) - \sum_{s=0}^{t-1} w^*_a(\hat{\bm \mu}(s))\right| \leq K (1+\sqrt{t}) \;.\]
\end{lemma}
It is slightly more efficient in practice to target directly $w^*(\hat{\bm \mu}(t))$, and to force exploration steps whenever an arm is in deficit. Introducing $U_t = \{ a : N_a(t) < \sqrt{t} - K/2\}$, the D-Tracking rule $(A_t)$ is sequentially defined as
\[A_{t+1} \in \left\{\begin{array}{ll}
                    \underset{a \in U_t}{\text{argmin}}  \ N_a(t) \ \text{if} \ U_t \neq \emptyset & (\textit{forced  exploration})\\
                    \underset{1 \leq a \leq K}{\text{argmax}} \ t \;w_a^*(\hat{\bm\mu}(t)) -{N_a}(t) & (\textit{direct tracking})
                   \end{array}
\right.\]

\begin{lemma}\label{lem:TrackingP} The D-Tracking rule ensures that $N_a(t) \geq (\sqrt{t} -K/2)_+ -1$ and that for all $\epsilon >0$, for all $t_0$, there exists $t_\epsilon \geq t_0$ such that 
\[\sup_{t\geq t_0}\max_a |w_a^*(\hat{\bm \mu}(t)) - w_a^*(\bm \mu)| \leq \epsilon \ \ \ \Rightarrow \ \ \ \sup_{t \geq t_\epsilon} \max_a \left|\frac{N_a(t)}{t} - w_a^*(\bm \mu)\right| \leq 3 (K-1)\epsilon\;.\]
\end{lemma}
It is guaranteed under all these sampling rules that the empirical proportion of draws of each arm converges to the optimal proportion, as proved in Appendix~\ref{proof:ConvergenceFraction}. 
\begin{proposition}\label{prop:ConvergenceFraction} The C-Tracking and D-Tracking sampling rules both satisfy
\[\bP_w\left(\lim_{t \rightarrow \infty} \frac{N_a(t)}{t} = w_a^*(\bm \mu)\right) = 1\;.\] 
\end{proposition}
We actually have a little more: a minimal convergence speed of $\hat{\mu}$ to $\mu$, which proves useful in the analysis of the~\emph{expected} sample complexity.
Of course, other tracking strategies are possible, like for example the one introduced in~\cite{CsabaTracking} for the uniform estimation of all arms' expectations. 

\subsection{Chernoff's Stopping Rule}\label{subsec:Stopping}
From a statistical point of view, the question of stopping at time $t$ is a more or less classical statistical test: do the past observations allow to assess, with a risk at most $\delta$, that one arm is larger than the others?
For all arms $a,b\in \Arms$, we consider the Generalized Likelihood Ratio statistic
\[Z_{a,b}(t) := \log\frac{\max_{\mu_a'\geq \mu_b'} p_{\mu_a'}\left(\underline{X}^a_{N_a(t)}\right)p_{\mu_b'}\left(\underline{X}^b_{N_b(t)}\right)}{\max_{\mu_a'\leq \mu_b'} p_{\mu_a'}\left(\underline{X}^a_{N_a(t)}\right)p_{\mu_b'}\left(\underline{X}^b_{N_b(t)}\right)}\;,\]
where $\underline{X}^a_{N_a(t)} = (X_s: A_s=a, s\leq t)$ is a vector that contains the observations of arm $a$ available at time $t$, and where $p_{\mu}(Z_1,\dots,Z_n)$ is the likelihood of $n$ i.i.d. observations from $w^{\mu}$: 
\[p_{\mu}(Z_1,\dots,Z_n) = \prod_{k=1}^{n} \exp(\dot{b}^{-1}(\mu) Z_k - b(\dot{b}^{-1}(\mu)))\;.\]
This statistic has a convenient closed-form expression for exponential family bandit models. Introducing for all arms $a,b$ a weighted average of their empirical mean: 
 \begin{eqnarray*}
\hat{\mu}_{a,b}(t) &:=& \frac{N_a(t)}{N_a(t) + N_b(t)}\hat{\mu}_a(t) + \frac{N_b(t)}{N_a(t) + N_b(t)}{\hat{\mu}_b(t)}\;,
\end{eqnarray*}
it is well known and easy to see that if $\hat{\mu}_a(t) \geq \hat{\mu}_b(t)$, 
	\begin{equation}Z_{a,b}(t) = N_a(t)\, d\big(\hat{\mu}_a(t), \hat{\mu}_{a,b}(t) \big) + N_b(t)\, d\big(\hat{\mu}_b(t), \hat{\mu}_{a,b}(t) \big)\;,\label{eq:Zformula}\end{equation}
and that $Z_{a,b}(t) = - Z_{b,a}(t)$.
The testing intuition thus suggests the following stopping rule: 
\begin{align}\tau_\delta &= \inf \left\{t \in \N : \exists a\in \Arms, \forall b \in \Arms \setminus \{a\}, Z_{a,b}(t)>\beta(t,\delta)\right\} \nonumber\\
&= \inf \left\{t \in \N : Z(t):= \max_{a\in \Arms} \min_{b \in \Arms \setminus \{a\}} Z_{a,b}(t)>\beta(t,\delta)\right\}\;,\label{def:Stopping}\end{align}
where $\beta(t,\delta)$ is an exploration rate to be tuned appropriately. The form of this stopping rule can be traced back to \cite{Chernoff59}\footnote{The stopping rule $\tau_\delta$ was proposed under equivalent form (with a different threshold) in the context of adaptive sequential hypothesis testing. Best arm identification in a bandit model can be viewed as a particular instance in which we test $K$ hypotheses, $H_a : ( \mu_a = \max_{i \in \Arms} \mu_i)$, based on adaptively sampling the marginal of $\nu = (\nu^{\mu_1},\dots,\nu^{\mu_K})$. However, \cite{Chernoff59} considers a different performance criterion, and its analysis holds when each of the hypotheses consists in a finite set of parameters, unlike the bandit setting.}.
As $\min_{b \in \Arms \setminus a} Z_{a,b}(t)$ is non-negative if and only if $\hat{\mu}_a(t)\geq \hat{\mu}_b(t)$ for all $b\neq a$, $Z(t) = \min_{b \in \Arms \setminus \hat{a}_t} Z_{\hat{a}_t,b}(t)$ whenever there is a unique  best empirical arm  $\hat{a}_t =  \text{argmax}_{a \in \Arms} \ \hat{\mu}_a(t)$ at time $t$.
 Obviously, $\big(\hat{\mu}_a(\tau_\delta)\big)_a$ has a unique maximizer, which is the final decision.

In addition to the testing interpretation given above, the stopping rule can be explained in light of the lower bound $\bE_\mu[\tau_\delta/T^*(\bm\mu)]\geq \kl(\delta,1-\delta)$. Indeed, one may write
\[Z(t) = t \min_{\lambda \in \mathrm{Alt}(\hat{\bm \mu}(t))}\sum_{a=1}^K \left(\frac{N_a(t)}{t}\right)d(\hat{\mu}_a(t),\lambda_a)  \leq \frac{t}{ T^*\big(\hat{\mu}(t)\big)}\;.\]
Theorem~\ref{cor:GeneBAI} suggests that a $\delta$-PAC strategy cannot stop (at least in expectation) before $t/T^*\big(\hat{\mu}(t)\big)$ is larger than $\kl(\delta, 1-\delta)\sim \log(1/\delta)$, which suggests to stop when $Z(t)\geq \log(1/\delta)$. In Section~\ref{sec:PACanalysis}, we prove that a slightly more prudent choice of the threshold $\beta(t,\delta)$ does lead to a PAC algorithm (whatever the sampling rule, even if the proportions of draws are sub-optimal). And it is shown in Section~\ref{sec:SCAnalysis} that, using our sampling rule, $\bE[\tau_\delta]$ is indeed of order $T^*(\bm\mu)\log(1/\delta)$.

It is also possible to give a Minimum Description Length (MDL) interpretation of the stopping rule. It is well known that choosing the model that gives the shortest description of the data is a provably efficient heuristic (see~\cite{Rissanen78}, and \cite{Grunwald07MDL} for a survey). In some sense, the stopping rule presented above follows the same principle. In fact, elementary algebra shows that 
\[Z_{a,b}(t) = (N_a(t) + N_b(t))h\left(\hat{\mu}_{a,b}(t)\right) - \left[N_a(t) h(\hat{\mu}_a(t)) + N_b(t)h(\hat{\mu}_b(t))\right],\]
where $h(\mu) = \bE_{X\sim \nu^\mu}[-\log p_\mu(X)] = b(\dot{b}^{-1}(\mu)) - \dot{b}^{-1}(\mu) \mu$. 
In the Bernoulli case, the \emph{Shannon entropy} $h(\mu)=-\mu\log(\mu)-(1-\mu)\log(1-\mu)$ is well-known to represent an ideal code length per character for binary compression. Thus, $Z_{a,b}(t)$ appears as the difference between the ideal code length for the rewards of arms $a$ and $b$ coded together, and the sum of the ideal code lengths for the rewards coded separately. If this difference is sufficiently large\footnote{Universal coding theory would suggest a threshold of order $\log(t)/2$, a term that will appear in Section~\ref{sec:PACanalysis}.}, the shortest description of the data is to separate arms $a$ and $b$. The stopping rule~\eqref{def:Stopping} thus consists in waiting for the moment when it becomes cheaper to code the rewards of the best empirical arm separately from each of the others.
It is no surprise that the proof of Proposition~\ref{prop:deltaPAC} below is based on a classical information-theoretic argument.

\section{Choosing the Threshold in the Stopping Rule}\label{sec:PACanalysis}
We now explain how to choose the threshold $\beta(t,\delta)$ so as to ensure the $\delta$-PAC property: with probability larger than $1-\delta$, any algorithm based on the stopping rule~\eqref{def:Stopping} outputs the optimal arm, provided that it stops. The interpretations of the stopping rule presented in the last section suggest the presence of two ingredients: $\log(1/\delta)$ for the risk, and $\log(t)$ for the fluctuations of the counts. 
We present here two results: one is based on an information-theoretic argument used for consistency proofs of MDL estimators, and the second is based on the probabilistic control of self-normalized averages taken from~\cite{Combes14Lip}. 
To keep things simple, the first argument is detailed only for the case of Bernoulli rewards (the standard framework of coding theory). 
The second argument is more general, but a little less tight.

\paragraph{The Informational Threshold.}
\begin{theorem}\label{prop:deltaPAC} 
Let $\delta \in (0,1)$. Whatever the sampling strategy, using Chernoff's stopping rule~\eqref{def:Stopping} on Bernoulli bandits with threshold 
\[\beta(t,\delta) = \log \left(\frac{2t (K-1)}{\delta}\right)\]
ensures that for all $\bm\mu \in \cS$, $\bP_{\bm\mu}\left(\tau_\delta<\infty, \,\hat{a}_{\tau_\delta} \neq a^*\right) \leq \delta.$
\end{theorem}
\paragraph{Proof sketch.} A more detailed proof is given in Appendix~\ref{sec:proof:deltaPAC}. We proceed here similarly to~\cite{Rajesh15Oddball}, employing an argument used for MDL consistency proofs (see~\cite{G06} and references therein).
Introducing, for $a,b\in\Arms$, $T_{a,b} := \inf\{ t\in \N : Z_{a,b}(t) > \beta(t,\delta)\}$, one has 
\begin{align*}
\bP_{\bm \mu}(\tau_\delta<\infty,\hat{a}_{\tau_\delta}\neq a^*) & \leq \bP_{\bm \mu}\left(\exists a \in \Arms\setminus a^*, \exists t \in \N : Z_{a,a^*}(t) > \beta(t,\delta)\right) \\ &\leq \sum_{a\in \Arms\setminus \{a^*\}} \bP_{\bm \mu}(T_{a,a^*} < \infty)\;.
\end{align*}
It is thus sufficient to show that if $\beta(t,\delta)=\log(2t/\delta)$, and if $\mu_a<\mu_b$, then $\bP_{\bm \mu}(T_{a,b} < \infty) \leq \delta$.
For such a pair of arms, observe that on the event $\big\{T_{a,b}=t\big\}$ time $t$ is the first moment when $Z_{a,b}$ exceeds the threshold $\beta(t,\delta)$, which implies by definition that  
\[1 \leq e^{-\beta(t,\delta)} \frac{\max_{{\mu}_a' \geq {\mu}'_b} p_{\mu_a'}(\underline{X}^a_t)p_{\mu_b'}(\underline{X}^b_t)}{\max_{{\mu}_a' \leq {\mu}'_b} p_{\mu_a'}(\underline{X}^a_t)p_{\mu_b'}(\underline{X}^b_t)}\;.\]
It thus holds that
 \begin{align*}
 &\bP_{\bm \mu}(T_{a,b}< \infty)  = \sum_{t=1}^\infty \bP_{\bm \mu}(T_{a,b}=t) = \sum_{t=1}^\infty \bE_{\bm \mu}\left[\ind_{(T_{a,b}=t)}\right] \\
 & \leq \sum_{t=1}^\infty e^{-\beta(t,\delta)}\bE_{\bm \mu}\left[\ind_{(T_{a,b}=t)} \frac{\max_{{\mu}_a' \geq {\mu}'_b} p_{\mu_a'}(\underline{X}^a_t)p_{\mu_b'}(\underline{X}^b_t)}{\max_{ {\mu}_a' \leq {\mu}'_b} p_{\mu_a'}(\underline{X}^a_t)p_{\mu_b'}(\underline{X}^b_t)}\right] \\
 & \leq \sum_{t=1}^\infty e^{-\beta(t,\delta)}\bE_{\bm \mu}\left[\ind_{(T_{a,b}=t)} \frac{\max_{{\mu}_a' \geq {\mu}'_b} p_{\mu_a'}(\underline{X}^a_t)p_{\mu_b'}(\underline{X}^b_t)}{p_{\mu_a}(\underline{X}^a_t)p_{\mu_b}(\underline{X}^b_t)}\right] \\
 & =  \sum_{t=1}^\infty e^{-\beta(t,\delta)}\int_{\{0,1\}^t}\ind_{(T_{a,b}=t)}(x_1,\dots,x_t) \underbrace{ \max_{{\mu}_a' \geq {\mu}'_b} p_{\mu_a'}(\underline{x}^a_t)p_{\mu_b'}(\underline{x}^b_t)\prod_{i \in \Arms \setminus\{a,b\}} p_{\mu_i}(\underline{x}_t^i)}_{(*)}dx_1\dots dx_t \;.
\end{align*}
Of course the maximum likelihood $(*)$ is not a probability density. A possible workaround (sometimes referred to as \emph{Barron's lemma}, see~\cite{BRY98} 
and references therein) is to use a \emph{universal distribution} like~\cite{KT81}, which is known to provide a tight uniform approximation:
\begin{lemma}\label{lem:ktMain}[\cite{Willems95thecontext}] Let $p_u(x)$ be the likelihood of successive observations $x \in \{0,1\}^n$ of a Bernoulli random variable with mean $u$. Then the Krichevsky-Trofimov distribution 
 \[\kt(x) = \int_0^1 \frac{1}{\pi\sqrt{u(1-u)}} p_u(x) \d u\;\]
is a probability law on $\{0,1\}^{n}$ that satisfies 
\[ \sup_{x\in\{0,1\}^n} \sup_{u \in[0,1]} \frac{p_u(x)}{\kt(x)}\leq 2\sqrt{n}\;.\]
 \end{lemma}
Together with the inequality $\sqrt{ab}\leq (a+b)/2$, this property permits to conclude that $\bP_{\bm \mu}(T_{a,b}< \infty)$ is upper-bounded by 
 \begin{align*}
    &\sum_{t=1}^\infty 2t e^{-\beta(t,\delta)}\int_{\{0,1\}^t}\ind_{(T_{a,b}=t)}(x_1,\dots,x_t) \overbrace{ \kt(\underline{x}^a_t)\kt(\underline{x}^b_t)\prod_{i \in \Arms \setminus\{a,b\}} p_{\mu_i}(\underline{x}_t^i)}^{:=I(x_1,\dots,x_t)}dx_1\dots dx_t \\
  &=  \sum_{t=1}^\infty 2t e^{-\beta(t,\delta)} \tilde{\bE}\left[\ind_{(T_{a,b}=t)}\right] \leq \delta \,\tilde{\bP}(T_{a,b} < \infty) \leq \delta,
\end{align*}
using that the partially integrated likelihood $I(x_1,\dots,x_t)$ is the density of a probability measure $\tilde{\bP}$ (we denote the corresponding expectation by $\tilde{\bE}$). \qed

\paragraph{The Deviational Threshold.} The universal coding argument above can be adapted to other distributions, as shown for example in~\cite{CGG09}. It is also possible to make use of a deviation result like~\cite{Combes14Lip} in order to prove PAC guarantees in any exponential family bandit model. The exploration rate involved are slightly larger and less explicit. 
\begin{proposition}\label{prop:PACGene} Let $\bm \mu$ be an exponential family bandit model. Let $\delta \in (0,1)$ and $\alpha>1$. There exists a constant $C=C(\alpha,K)$ such that whatever the sampling strategy, using Chernoff's stopping rule~\eqref{def:Stopping} with the threshold
\[\beta(t,\delta)=\log \left(\frac{C t^\alpha}{\delta}\right)\]
ensures that for all $\bm\mu \in \cS$, $\bP_{\bm\mu}\big(\tau_\delta<\infty, \, \hat{a}_{\tau_\delta} \neq a^*\big) \leq \delta.$
\end{proposition}
The proof of Proposition~\ref{prop:PACGene} is given in Appendix~\ref{proof:PACGene}. 

\section{Sample Complexity Analysis}\label{sec:SCAnalysis}
Combining Chernoff's stopping rule and an optimal-tracking sampling rule permits to approach the lower bound of Theorem~\ref{cor:GeneBAI} for sufficiently small values of the risk $\delta$. 
We first state a simple almost-sure convergence result and give its short proof. 
Then, we sketch the (somewhat more technical) analysis controlling  the expectation of $\tau_\delta$.
\subsection{Almost-sure Upper Bound}
\begin{proposition}\label{lem:AnalysisAS} Let $\alpha\in [1,e/2]$ and $r(t) = O(t^\alpha)$. Using Chernoff's stopping rule with $\beta(t,\delta)=\log(r(t)/\delta)$, and any sampling rule ensuring that for every arm $a\in\Arms$, $N_a(t)/t$ converges almost-surely to $w_a^*$ guarantees that for all $\delta\in (0,1)$, $\bP_{\bm\mu}(\tau_\delta < +\infty)=1$ and  \[\bP_{\bm\mu}\left(\limsup_{\delta \rightarrow 0}\frac{\tau_\delta}{\log(1/\delta)} \leq \alpha T^*(\bm \mu) \right) = 1.\]
 \end{proposition}
\paragraph{Proof.}
Let $\cE$ be the event 
\[\cE = \left\{\forall a \in \Arms, \frac{N_a(t)}{t} \underset{t \rightarrow \infty}{\rightarrow} w_a^* \ \ \text{and} \ \ \ \hat{\bm \mu}(t) \underset{t \rightarrow \infty}{\rightarrow} \bm \mu\right\}\;.\]
From the assumption on the sampling strategy and the law of large number, $\cE$ is of probability 1.
On $\cE$, there exists $t_0$ such that for all $t \geq t_0$, $\hat{\mu}_1(t) > \max_{a \neq 1} \hat{\mu}_a(t)$ and
\begin{eqnarray*}
 Z(t) & = & \min_{a \neq 1} \ Z_{1,a}(t) =  \min_{a \neq 1} N_1(t)d(\hat{\mu}_1(t) , \hat{\mu}_{1,a}(t)) + N_a(t) d(\hat{\mu}_a(t),\hat{\mu}_{1,a}(t))\\
 & = & t \left[\min_{a \neq 1} \left(\frac{N_1(t)}{t} + \frac{N_a(t)}{t}\right)I_{\frac{N_1(t)/t}{N_1(t)/t + N_a(t)/t}}\left(\hat{\mu}_1(t),\hat{\mu}_a(t)\right)  \right]. 
\end{eqnarray*}
For all $a \geq 2$, the mapping $(\bm w , \bm \lambda) \rightarrow (w_1 + w_a) I_{w_1/(w_1 + w_a)}(\lambda_1,\lambda_a)$ is continuous at $(w^*(\bm \mu),\bm \mu)$. 
Therefore, for all $\epsilon > 0$ there exists $t_1 \geq t_0$ such that for all $t \geq t_1$ and all $a\in \{2,\dots,K\}$,
\[\left(\frac{N_1(t)}{t} + \frac{N_a(t)}{t}\right)I_{\frac{N_1(t)/t}{N_1(t)/t + N_a(t)/t}}\left(\hat{\mu}_1(t),\hat{\mu}_a(t)\right) \geq \frac{w^*_1 + w^*_a}{1+\epsilon}I_{\frac{w_1^*}{w_1^*+ w_a^*}}\left({\mu}_1,{\mu}_a\right).\]
Hence, for $t \geq t_1$, 
\[Z(t) \geq \frac{t}{1+\epsilon} \min_{a \neq 1}\left(w^*_1 + w^*_a\right)I_{\frac{w_1^*}{w_1^*+ w_a^*}}\left({\mu}_1,{\mu}_a\right) = \frac{t}{(1+\epsilon)T^*(\bm \mu)}.\]
Consequently,
\begin{eqnarray*}
 \tau_\delta &= &  \inf\{ t \in \N : Z(t) \geq \beta(t,\delta) \} \\
 & \leq & t_1 \vee \inf\{t \in \N : t (1+\epsilon)^{-1}T^*(\bm \mu)^{-1} \geq \log(r(t) /\delta) \} \\ 
 & \leq & t_1 \vee \inf\{t \in \N : t (1+\epsilon)^{-1}T^*(\bm \mu)^{-1} \geq \log(Ct^{\alpha} /\delta) \}, 
\end{eqnarray*}
for some positive constant $C$. Using the technical Lemma~\ref{lem:technical} in the Appendix, it follows that on $\cE$, as $\alpha \in [1,e/2]$,
\[\tau_\delta \leq t_1 \vee \alpha (1+\epsilon) T^*(\bm \mu) \left[\log \left(\frac{Ce((1+\epsilon)T^*(\bm \mu))^\alpha}{\delta}\right) +  \log\log \left(\frac{C((1+\epsilon)T^*(\bm \mu))^\alpha}{\delta}\right)\right]\;.\]
Thus $\tau_\delta$ is finite on $\cE$ for every $\delta \in (0,1)$, and 
\[\limsup_{\delta \rightarrow 0} \frac{\tau_\delta}{\log(1/\delta)} \leq (1+\epsilon)\,\alpha\,  T^*(\bm \mu)\;.\]
Letting $\epsilon$ go to zero concludes the proof.
\qed

\subsection{Asymptotic Optimality in Expectation}
In order to prove that the lower bound of Theorem~\ref{cor:GeneBAI} is matched, we now give an upper bound on the expectation of the stopping time $\tau_\delta$. The proof of this result is to be found in Appendix~\ref{sec:ProofSC}.

\begin{theorem}\label{thm:AsymptoticSC} Let $\bm\mu$ be an exponential family bandit model. Let $\alpha\in [1,e/2]$ and $r(t) = O(t^\alpha)$. Using Chernoff's stopping rule with $\beta(t,\delta)=\log(r(t)/\delta)$, and the sampling rule C-Tracking or D-Tracking,
\[\limsup_{\delta \rightarrow 0} \frac{\bE_{\bm\mu}[\tau_\delta]}{\log(1/\delta)} \leq {\alpha}{T^*(\bm\mu)}\;.\] 
\end{theorem}
To summarize, for every Bernoulli bandit $\bm\mu$, the choice $\beta(t,\delta)=\log(2(K-1)t/\delta)$ in Chernoff's stopping rule is $\delta$-PAC (by Theorem~\ref{prop:deltaPAC}); with one of the sampling rules given above, the stopping time $\tau_\delta$ is almost surely finite (by Proposition~\ref{lem:AnalysisAS}) and when $\delta$ is small enough its expectation is close to $T^*(\bm\mu)\log(1/\delta)$ by Theorem~\ref{thm:AsymptoticSC}, an optimal sample complexity after Theorem~\ref{cor:GeneBAI}. 

More generally, for exponential family bandit models, combining Proposition~\ref{prop:PACGene} and Theorem~\ref{thm:AsymptoticSC}, one obtains for every $\alpha>1$ the existence of an exploration rate for which Chernoff's stopping rule combined with a tracking sampling rule is  $\delta$-PAC and satisfies
\[\limsup_{\delta \rightarrow 0} \frac{\bE_{\bm\mu}[\tau_\delta]}{\log(1/\delta)} \leq \alpha{T^*(\bm\mu)}\;.\] 

\section{Discussion and Numerical Experiments}\label{sec:experiments}
We give here a few comments on the practical behaviour of the Track-and-Stop (T-a-S) strategy.
Let us first emphasize that the forced exploration step are rarely useful, but in some cases really necessary and not only for the theorems: when $\mu_2$ and $\mu_3$ are equal, they prevent the probability that the strategy never ends from being strictly positive. 
Second, our simulation study suggests that the exploration rate $\beta(t,\delta)=\log((\log(t)+1)/\delta)$, though not (yet) allowed by theory, is still over-conservative in practice. Further, even though any sampling strategy ensuring that $N_a(t) / t \to w_a^*$  satisfies the optimality theorems above, we propose (without formal justification) an experimentally more efficient sampling rule: after $t$ observations, let 
\[\hat{c}_t = \argmin{c\in\Arms\setminus\{\hat{a}_t\}} \,  Z_{\hat{a}_t,c}(t) \] 
be the  'best challenger' of the empirical champion $\hat{a}_t$. 
 We choose $A_{t+1} = \hat{a}_t$~if 
\[\frac{N_{\hat{a}_t}(t)}{N_{\hat{a}_t}(t) + N_{\hat{c}_t}(t)} < \frac{w_{\hat{a}_t}^*(\hat{\bm \mu}(t))}{w_{\hat{a}_t}^*(\hat{\bm \mu}(t))+w_{\hat{c}(t)}^*(\hat{\bm \mu}(t))}\]
and $A_{t+1} = \hat{c}_t$ otherwise (with forced explorations steps as in the D-Tracking rule).

We consider two sample scenarios $\bm\mu_1=[0.5 \ 0.45 \ 0.43 \ 0.4]$ and $\bm\mu_2=[0.3 \ 0.21 \ 0.2 \ 0.19 \ 0.18]$, and we choose $\delta=0.1$. This choice is meant to illustrate that our algorithm performs well even for relatively high risk values (so far, optimality is proved only for small risks).
We compare the Track-and-Stop algorithm based on D-Tracking and BestChallenger to algorithms from the literature designed for Bernoulli bandit models, namely KL-LUCB and KL-Racing~\citep{COLT13}. Racing algorithms proceed in rounds: at start, all arms are active; at each round, all active arms are drawn once; at the end of a round, a rule determines if the empirically worst arm should be eliminated. Call $\hat{a}_r$ the empirically best arm after $r$ rounds. In KL-Racing, arm $b$ is eliminated if its upper confidence bound $ \max \{q \in [0,1] : r d(\hat{\mu}_{b,r},q ) \leq \beta(r,\delta)\}$ is smaller than the best arm's lower bound $\min \{q \in [0,1] : r d(\hat{\mu}_{\hat{a}_r,r},q ) \leq \beta(r,\delta)\}$. We also introduced in the competition the `hybrid' Chernoff-Racing strategy, which eliminates $b$ if  \[Z_{\hat{a}_r, b} = rd\left(\hat{\mu}_{\hat{a}_r,r} , \frac{\hat{\mu}_{\hat{a}_r,r}  + \hat{\mu}_{b,r} }{2}\right) + rd\left(\hat{\mu}_{b,r} , \frac{\hat{\mu}_{\hat{a}_r,r}  + \hat{\mu}_{b,r} }{2}\right) > \beta(r,\delta).\]
Table~\ref{tab:results} presents the estimated average number of draws of the five algorithms in the two scenarios. Our (Julia) code will be available online. We see that the use of the MDL stopping rule leads to a clear improvement. Moreover, Chernoff-Racing significantly improves over KL-Racing, and its performance is even close to that of our optimal algorithms. 

\begin{table}[h]
\centering
\begin{tabular}{|c|c|c|c|c|c|}
\hline
 & T-a-S (BC) & T-a-S (D-Tracking) &  Chernoff-Racing & KL-LUCB & KL-Racing \\
\hline
$\bm\mu_1$ 
&  3968  & 4052   & 4516 & 8437 &  9590 \\\hline
$\bm \mu_2$
& 1370 & 1406 & 3078   & 2716 &  3334 \\
\hline
\end{tabular}
\caption{\label{tab:results} Expected number of draws $\bE_{\bm\mu}[\tau_\delta]$ for $\delta=0.1$, averaged over $N=3000$ experiments:
$\bm\mu_1=[0.5 \ 0.45 \ 0.43 \ 0.4]$, $w^*(\bm \mu_1) =[0.417 \ 0.390 \ 0.136 \ 0.057]$;
$\bm\mu_2=[0.3 \ 0.21 \ 0.2 \ 0.19 \ 0.18]$, $w^*(\bm\mu_2)=[0.336 \ 0.251 \ 0.177  \ 0.132 \ 0.104]$.
}
\end{table}
It should be emphasized that a Racing-type algorithm cannot reach the lower bound in general: by construction, it forces the last two arms in the race (hopefully $\mu_1$ and $\mu_2$) to be drawn equally often, which is sub-optimal unless $w_1^*(\bm \mu) = w^*(\bm \mu_2)$ (a condition approximately matched only if there is a large gap between the second and the third best arms). This is illustrated in the second scenario $\bm\mu_2$ of Table~\ref{tab:results}, where the ratio $w_1^*/w_2^*$ is larger.

\section{Conclusion}
We gave a characterization of the complexity of best arm identification in the fixed confidence-setting, for a large class of bandit models with arms parameterized by their means. Our new lower bound reveals the existence of optimal proportions of draws of the arms that can be computed efficiently. Our  Track-and-Stop strategy, by combining a track of these optimal proportions with Chernoff's stopping rule, asymptotically matches the lower bound. In future work, instead of using forced exploration steps within a plugin procedure, we will investigate optimistic (or robust-to-noise) sampling strategies in order to optimize the exploration and to obtain non-asymptotic sample complexity bounds. Furthermore, we will investigate the fixed-budget setting, for which we conjecture that $P(\hat{a}_t \neq 1) \gtrsim \exp\big(-t/T_*(\bm\mu)\big)$ with
\[T_*(\bm\mu)^{-1}  = \sup_{w\in\Sigma_K} \min_{a\in\{2,\dots,K\}} \inf_{\mu_a<m_a<\mu_1} w_1\, d(m_a, \mu_1) + w_a\, d(m_a, \mu_a)\;.\]

\acks{
This work was partially supported by the CIMI (Centre International de Math\'ematiques et d'{In\-for\-ma\-tique}) Excellence program while Emilie Kaufmann visited Toulouse in November 2015.
The authors acknowledge the support of the French Agence Nationale de la Recherche (ANR), under grants ANR-13-BS01-0005 (project SPADRO) and ANR-13-CORD-0020 (project ALICIA). They thank S\'ebastien Gerchinovitz for stimulating discussions and comments.
}

\bibliography{biblioBandits}

\begin{thebibliography}{29}
\providecommand{\natexlab}[1]{#1}
\providecommand{\url}[1]{\texttt{#1}}
\expandafter\ifx\csname urlstyle\endcsname\relax
  \providecommand{\doi}[1]{doi: #1}\else
  \providecommand{\doi}{doi: \begingroup \urlstyle{rm}\Url}\fi

\bibitem[Abbasi-Yadkori et~al.(2011)Abbasi-Yadkori, D.P{\'a}l, and
  C.Szepesv{\'a}ri]{YadLinear11}
Y.~Abbasi-Yadkori, D.P{\'a}l, and C.Szepesv{\'a}ri.
\newblock {Improved Algorithms for Linear Stochastic Bandits}.
\newblock In \emph{{Advances in Neural Information Processing Systems}}, 2011.

\bibitem[Agrawal and Goyal(2013)]{AGContext13}
S.~Agrawal and N.~Goyal.
\newblock {Thompson Sampling for Contextual Bandits with Linear Payoffs}.
\newblock In \emph{{International Conference on Machine Learning (ICML)}},
  2013.

\bibitem[Antos et~al.(2008)Antos, Grover, and Szepesv\'{a}ri]{CsabaTracking}
A.~Antos, V.~Grover, and C.~Szepesv\'{a}ri.
\newblock Active learning in multi-armed bandits.
\newblock In \emph{Algorithmic Learning Theory}, 2008.

\bibitem[Barron et~al.(1998)Barron, Rissanen, and Yu]{BRY98}
A.~Barron, J.~Rissanen, and Bin Yu.
\newblock The minimum description length principle in coding and modeling.
\newblock \emph{Information Theory, IEEE Transactions on}, 44\penalty0
  (6):\penalty0 2743--2760, Oct 1998.
\newblock ISSN 0018-9448.
\newblock \doi{10.1109/18.720554}.

\bibitem[Bubeck and Cesa-Bianchi(2012)]{Bubeck:Survey12}
S.~Bubeck and N.~Cesa-Bianchi.
\newblock {Regret analysis of stochastic and nonstochastic multi-armed bandit
  problems}.
\newblock \emph{Fondations and Trends in Machine Learning}, 5(1):\penalty0
  1--122, 2012.

\bibitem[Bubeck et~al.(2011)Bubeck, Munos, Stoltz, and
  Szepesv\'{a}ri]{Bubeck11Xarmed}
S.~Bubeck, R.~Munos, G.~Stoltz, and C.~Szepesv\'{a}ri.
\newblock X-armed bandits.
\newblock \emph{Journal of Machine Learning Research}, 12:\penalty0 1587--1627,
  2011.

\bibitem[Capp{\'e} et~al.(2013)Capp{\'e}, Garivier, Maillard, Munos, and
  Stoltz]{KLUCBJournal}
O.~Capp{\'e}, A.~Garivier, O-A. Maillard, R.~Munos, and G.~Stoltz.
\newblock {{K}ullback-{L}eibler upper confidence bounds for optimal sequential
  allocation}.
\newblock \emph{Annals of Statistics}, 41(3):\penalty0 1516--1541, 2013.

\bibitem[Chambaz et~al.(2009)Chambaz, Garivier, and Gassiat]{CGG09}
Antoine Chambaz, Aur\'elien Garivier, and Elisabeth Gassiat.
\newblock A {MDL} approach to {HMM} with poisson and gaussian emissions.
  application to order identification.
\newblock \emph{Journal of Statistical Planning and Inference}, 139\penalty0
  (3):\penalty0 962--977, 2009.

\bibitem[Chernoff(1959)]{Chernoff59}
H.~Chernoff.
\newblock {Sequential design of Experiments}.
\newblock \emph{The Annals of Mathematical Statistics}, 30\penalty0
  (3):\penalty0 755--770, 1959.

\bibitem[Combes and Prouti{\`e}re(2014)]{Combes14Unimodal}
R.~Combes and A.~Prouti{\`e}re.
\newblock {Unimodal Bandits without Smoothness}.
\newblock Technical report, 2014.

\bibitem[Even-Dar et~al.(2006)Even-Dar, Mannor, and Mansour]{EvenDaral06}
E.~Even-Dar, S.~Mannor, and Y.~Mansour.
\newblock {Action Elimination and Stopping Conditions for the Multi-Armed
  Bandit and Reinforcement Learning Problems}.
\newblock \emph{Journal of Machine Learning Research}, 7:\penalty0 1079--1105,
  2006.

\bibitem[Gabillon et~al.(2012)Gabillon, Ghavamzadeh, and
  Lazaric]{Gabillon:al12}
V.~Gabillon, M.~Ghavamzadeh, and A.~Lazaric.
\newblock {Best Arm Identification: A Unified Approach to Fixed Budget and
  Fixed Confidence}.
\newblock In \emph{{Advances in Neural Information Processing Systems}}, 2012.

\bibitem[Garivier(2006)]{G06}
Aurélien Garivier.
\newblock Consistency of the unlimited {BIC} context tree estimator.
\newblock \emph{IEEE Transactions on Information Theory}, 52\penalty0
  (10):\penalty0 4630--4635, 2006.

\bibitem[Graves and Lai(1997)]{GravesLai97}
T.L. Graves and T.L. Lai.
\newblock {Asymptotically Efficient adaptive choice of control laws in
  controlled markov chains}.
\newblock \emph{SIAM Journal on Control and Optimization}, 35(3):\penalty0
  715--743, 1997.

\bibitem[Gr\"{u}nwald(2007)]{Grunwald07MDL}
Peter~D. Gr\"{u}nwald.
\newblock \emph{The Minimum Description Length Principle (Adaptive Computation
  and Machine Learning)}.
\newblock The MIT Press, 2007.
\newblock ISBN 0262072815.

\bibitem[Jamieson et~al.(2014)Jamieson, Malloy, Nowak, and
  Bubeck]{Jamiesonal14LILUCB}
K.~Jamieson, M.~Malloy, R.~Nowak, and S.~Bubeck.
\newblock {lil'{UCB}: an Optimal Exploration Algorithm for Multi-Armed
  Bandits}.
\newblock In \emph{{Proceedings of the 27th Conference on Learning Theory}},
  2014.

\bibitem[Kalyanakrishnan et~al.(2012)Kalyanakrishnan, Tewari, Auer, and
  Stone]{Shivaramal12}
S.~Kalyanakrishnan, A.~Tewari, P.~Auer, and P.~Stone.
\newblock {{PAC} subset selection in stochastic multi-armed bandits}.
\newblock In \emph{{International Conference on Machine Learning (ICML)}},
  2012.

\bibitem[Kaufmann and Kalyanakrishnan(2013)]{COLT13}
E.~Kaufmann and S.~Kalyanakrishnan.
\newblock {Information complexity in bandit subset selection}.
\newblock In \emph{{Proceeding of the 26th Conference On Learning Theory.}},
  2013.

\bibitem[Kaufmann et~al.(2014)Kaufmann, Capp{\'e}, and Garivier]{COLT14}
E.~Kaufmann, O.~Capp{\'e}, and A.~Garivier.
\newblock {On the Complexity of A/B Testing}.
\newblock In \emph{{Proceedings of the 27th Conference On Learning Theory}},
  2014.

\bibitem[Kaufmann et~al.(2015)Kaufmann, Capp{\'e}, and Garivier]{JMLR15}
E.~Kaufmann, O.~Capp{\'e}, and A.~Garivier.
\newblock {On the Complexity of Best Arm Identification in Multi-Armed Bandit
  Models}.
\newblock \emph{Journal of Machine Learning Research (to appear)}, 2015.

\bibitem[Krichevsky and Trofimov(1981)]{KT81}
Raphail~E. Krichevsky and Victor~K. Trofimov.
\newblock The performance of universal encoding.
\newblock \emph{{IEEE} Transactions on Information Theory}, 27\penalty0
  (2):\penalty0 199--206, 1981.
\newblock \doi{10.1109/TIT.1981.1056331}.
\newblock URL \url{http://dx.doi.org/10.1109/TIT.1981.1056331}.

\bibitem[Lai and Robbins(1985)]{LaiRobbins85bandits}
T.L. Lai and H.~Robbins.
\newblock {Asymptotically efficient adaptive allocation rules}.
\newblock \emph{Advances in Applied Mathematics}, 6\penalty0 (1):\penalty0
  4--22, 1985.

\bibitem[Magureanu et~al.(2014)Magureanu, Combes, and
  Prouti{\`e}re]{Combes14Lip}
S.~Magureanu, R.~Combes, and A.~Prouti{\`e}re.
\newblock {Lipschitz Bandits: Regret lower bounds and optimal algorithms}.
\newblock In \emph{{Proceedings on the 27th Conference On Learning Theory}},
  2014.

\bibitem[Mannor and Tsitsiklis(2004)]{MannorTsi04}
S.~Mannor and J.~Tsitsiklis.
\newblock {The Sample Complexity of Exploration in the Multi-Armed Bandit
  Problem}.
\newblock \emph{Journal of Machine Learning Research}, pages 623--648, 2004.

\bibitem[Munos(2014)]{SurveyRemiMCTS}
R.~Munos.
\newblock \emph{From bandits to Monte-Carlo Tree Search: The optimistic
  principle applied to optimization and planning.}, volume~7.
\newblock Foundations and Trends in Machine Learning, 2014.

\bibitem[Rissanen(1978)]{Rissanen78}
J.~Rissanen.
\newblock Modeling by shortest data description.
\newblock \emph{Automatica}, 14\penalty0 (5):\penalty0 465 -- 471, 1978.
\newblock ISSN 0005-1098.
\newblock \doi{http://dx.doi.org/10.1016/0005-1098(78)90005-5}.
\newblock URL
  \url{http://www.sciencedirect.com/science/article/pii/0005109878900055}.

\bibitem[Srinivas et~al.(2010)Srinivas, Krause, Kakade, and
  Seeger]{SrinivasGPUCB}
N.~Srinivas, A.~Krause, S.~Kakade, and M.~Seeger.
\newblock {Gaussian Process Optimization in the Bandit Setting : No Regret and
  Experimental Design}.
\newblock In \emph{{Proceedings of the International Conference on Machine
  Learning}}, 2010.

\bibitem[Vaidhyan and Sundaresan(2015)]{Rajesh15Oddball}
N.K. Vaidhyan and R.~Sundaresan.
\newblock Learning to detect an oddball target.
\newblock \emph{arXiv:1508.05572}, 2015.

\bibitem[Willems et~al.(1995)Willems, Shtarkov, and
  Tjalkens]{Willems95thecontext}
Frans M.~J. Willems, Yuri~M. Shtarkov, and Tjalling~J. Tjalkens.
\newblock The context tree weighting method: Basic properties.
\newblock \emph{IEEE Transactions on Information Theory}, 41:\penalty0
  653--664, 1995.

\end{thebibliography}

\appendix

\section{Characterization of the optimal proportion of draws}\label{proofs:Nu}

\subsection{Proof of Lemma~\ref{lem:CDinterpretation}}

Let $\bm \mu$ such that $\mu_1>\mu_2 \geq \dots \geq \mu_K$.  Using the fact that 
\[\mathrm{Alt}(\bm\mu)=\bigcup_{a \neq 1} \left\{\bm \lambda \in \cS : \lambda_a > \lambda_{1}\right\},\]
one has 
\begin{eqnarray*}
T^*(\bm\mu)^{-1} & = & \sup_{w \in \Sigma_K} \min_{a \neq 1} \inf_{\bm \lambda \in \cS : \lambda_a > \lambda_{1}} \sum_{a=1}^K w_a d(\mu_a,\lambda_a) \\
& = &  \sup_{w \in \Sigma_K} \min_{a \neq 1} \inf_{\bm \lambda  : \lambda_a \geq \lambda_{1}} \left[w_1d(\mu_{1},\lambda_{1})+w_a d(\mu_a,\lambda_a)\right]\;.
\end{eqnarray*}
Minimizing
\[f(\lambda_{1},\lambda_a) = w_1d(\mu_{1},\lambda_{1})+w_a d(\mu_a,\lambda_a)\]
under the constraint $\lambda_a \geq \lambda_{1}$ is a convex optimization problem that can be solved analytically. The minimum is obtained for 
\[\lambda_{1} = \lambda_a = \frac{w_1}{w_1+w_a}\mu_{1} + \frac{w_a}{w_1+w_a}\mu_{a}\]
and its value can be rewritten $(w_1 + w_a) I_{\frac{w_1}{w_1 + w_a}}(\mu_{1},\mu_a)$, using the function $I_\alpha$ defined in~\eqref{def:Divergences}.

\subsection{Proof of Theorem~\ref{thm:ExplicitForm}}\label{subsec:proofExplicit}
The function $g_a$ introduced in~\eqref{def:InverseFunction} rewrites
\[g_a(x)  =  d\left(\mu_1, m_a(x)\right) + x d(\mu_a,m_a(x)), \ \ \text{with} \ \ m_a(x) = \frac{\mu_1 + x \mu_a}{1+x}\;.\]
Using that $m_a'(x) = (\mu_a - \mu_1)/(1+x)^2$ and $\frac{d}{dy}d(x,y) = (y-x)/\ddot{b}(b^{-1}(y))$ one can show that $g_a$ is strictly increasing, since $g_a'(x) = d(\mu_a,m_a(x))>0$.
As $g_a(x)$ tends to $d(\mu_1,\mu_a)$ when $x$ goes to infinity, the inverse function $x_a(y) = g_a^{-1}(y)$ is defined on $[0,d(\mu_1,\mu_a)[$ and satisfies 
\[x_a'(y) = \frac{1}{d(\mu_a,m_a(x_a(y)))}>0\;.\]
Let $\bm w^*$ be an element in 
\begin{eqnarray*}
 \argmax{w \in \Sigma_K} \ \min_{a \neq 1} \ (w_1 + w_a)I_{\frac{w_1}{w_1+w_a}}(\mu_{1},\mu_a) 
 & = &\argmax{w \in \Sigma_K} \ w_1 \min_{a \neq 1} g_a\left(\frac{w_a}{w_1}\right)\;.
\end{eqnarray*}
The equality uses that $w_1^*\neq 0$ (since if $w_1=0$, the value of the objective is zero). Introducing $x_a^* = \frac{w_a^*}{w_1^*}$ for all $a\neq 1$, one has 
\[w_1^* = \frac{1}{1 + \sum_{a=2}^K x_a^*} \ \ \ \ \text{and,  for} \ a\geq 2, \ \ w_a^* =  \frac{x_a^*}{1 + \sum_{a=2}^K x_a^*}\]
and $(x_2^*,\dots,x_K^*) \in \R^{K-1}$ belongs to 
\begin{eqnarray}
\argmax{(x_2,\dots,x_K) \in \R^{K-1}} \ \frac{\min_{a \neq 1} g_a\left(x_a\right)}{1 + x_2 + \dots + x_{K}}. \label{equ:Argmax}
\end{eqnarray}
We now show that all the $g_a(x_a^*)$ have to be equal (Lemma~\ref{lem:allEqual}). Let \[\cB = \left\{ b \in \{2,\dots,K\} : g_b(x_b^*) = \min_{a \neq 1} g_a(x_a^*)\right\}\] and $\cA = \{2,\dots,K\} \backslash \cB$. Assume that $\cA \neq \emptyset$. For all $a\in \cA$ and $b\in \cB$, one has $g_a(x_a^*) > g_b(x_b^*)$. Using the continuity of the $g$ functions and the fact that they are strictly increasing, there exists $\epsilon>0$ such that 
\[\forall a \in \cA, b \in \cB, \ \ \  g_a(x_a^*-{\epsilon}/|\cA|)) > g_b\left(x_b^* + {\epsilon}/{|\cB|}\right)>g_b\left(x_b^*\right).\]
Introducing $\overline{x}_a = x_a^*-\epsilon/|\cA|$ for all $a\in\cA$ and $\overline{x}_b =x_b^* + {\epsilon}/{|\cB|}$ for all $b \in \cB$, there exists $b\in \cB$: 
\[\frac{\min_{a \neq 1} g_a\left(\overline{x}_a\right)}{1 + \overline{x}_2 + \dots \overline{x}_K}  = \frac{g_b\left(x_b^* + {\epsilon}/{|\cB|}\right)}{1 + x_2^* + \dots + x^*_{K}}  >  \frac{g_b\left(x_b^*\right)}{1 + x_2^* + \dots + x^*_{K}} =\frac{\min_{a\neq 1} g_a(x_a^*)}{1 + x_2^* + \dots + x^*_{K}}\;,\]
which contradicts the fact that $\bm x^*$ belongs to \eqref{equ:Argmax}. Hence $\cA=\emptyset$ and there exists $y^*\in [0,d(\mu_1,\mu_2)[$ such that 
\[\forall a \in \{2,\dots,K\}, \ g_a(x_a^*) = y^* \ \ \Leftrightarrow \ \ x_a^* = x_a(y^*),\]
with the function $x_a$ introduced above. From \eqref{equ:Argmax}, $y^*$ belongs to 
\begin{eqnarray*}
\argmax{y \in [0,d(\mu_1,\mu_2)[} \ G(y) \ \ \ \text{with} \ \ \ G(y)=\frac{y}{1 + x_2(y) + \dots + x_{K}(y)}\;. 
\end{eqnarray*}
$G$ is differentiable and, using the derivative of the $x_a$ given above, $G'(y)=0$ is equivalent to 
\begin{eqnarray}
 \sum_{a=2}^K\frac{y}{d(\mu_a,m_a(x_a(y))}&= &  1 + x_2(y) + \dots + x_{K}(y)\nonumber\\
\sum_{a=2}^K\frac{d(\mu_1,m_a(x_a(y))) + x_a(y)d(\mu_a,m_a(x_a(y))}{d(\mu_a,m_a(x_a(y))}  &= & 1 + x_2(y) + \dots + x_{K}(y)\nonumber\\
 \sum_{a=2}^K\frac{d(\mu_1,m_a(x_a(y)))}{d(\mu_a,m_a(x_a(y))} & = & 1\;.\label{equ:SolY}
\end{eqnarray}
For the the second equality, we use that $\forall a ,  d(\mu_1,m_a(x_a(y))) + x_a(y)d(\mu_a,m_a(x_a(y)) = y$. Thus $y^*$ is solution of the equation \eqref{equ:SolY}. This equation has a unique solution since 
\begin{equation}F_{\bm \mu}(y) = \sum_{a=2}^K \frac{d(\mu_1,m_a(x_a(y)))}{d(\mu_a,m_a(x_a(y))}\label{functionF}\end{equation}
is strictly increasing and satisfies $F_{\bm \mu}(0)=0$ and $\lim_{y \rightarrow d(\mu_1,\mu_2)} F_{\bm \mu}(y) = +\infty$. As $G$ is positive and satisfies $G(0)=0$, $\lim_{y \rightarrow d(\mu_1,\mu_2)} G(y)=0$, the unique local extrema obtained in $y^*$ is a maximum. 

\subsection{Proof of Proposition~\ref{prop:PropNu}}

Let $\bm \mu \in \cS$, and re-label its arms in decreasing order. From Theorem~\ref{thm:ExplicitForm}, $w_{1}(\bm \mu) \neq 0$ and for $a\geq 2$, $w_a^*(\bm \mu) =0 \Leftrightarrow x_a(y^*)=0 \Rightarrow y^*=0$ where $y^*$ is the solution of \eqref{equ:SolY}. But $0$ is not solution of \eqref{equ:SolY}, since the value of the left-hand side is 0. This proves that for all $a$, $w_a^*(\bm \mu) \neq 0$. 
For a given $\bm\mu$, $y^*$ is defined by   
\[F_{\bm\mu}(y^*) - 1 =0\;,\]
where $F_{\bm\mu}$ is defined in \eqref{functionF}. For all $\bm \mu \in \cS$ and every $y \in [0, d(\mu_1,\mu_2)[$, it can be shown that $\frac{d}{dy}F_{\bm\mu}(y) \neq 0$, in particular $\frac{d}{dy}F_{\bm\mu}(y^*)\neq 0$. Thus $y^*$ is a function of $\bm\mu$ that is continuous in every $\bm\mu \in \cS$, denoted by $y^*(\bm\mu)$. By composition, the function $\bm \mu \mapsto x_a(y^*(\bm\mu))$ are continuous in $\bm \mu \in  \cS$, and so does $w^*$.    

The proof of Statement 3 relies on the fact that if $a$ and $b$ are such that $\mu_1 > \mu_a \geq \mu_b$, $g_a(x) \leq g_b(x)$ for all $x$. Thus, for all $y \in [0,d(\mu_1,\mu_2)[$, $x_a(y) \geq x_b(y)$ and particularizing for $y^*(\bm \mu)$ yields the result. 

\subsection{Bounds on the characteristic time in the Gaussian case}\label{proof:BoundGaussian}

In the Gaussian case, with $d(x,y)=(x-y)^2/(2\sigma^2)$, the expression in Lemma~\ref{lem:CDinterpretation} can be made more explicit and yields 
\[T^*(\bm\mu)^{-1} = \sup_{w \in \Sigma_K}\inf_{a \neq 1} \frac{w_1w_a}{w_1 + w_a}\frac{\Delta_a^2}{2\sigma^2}\;.\]
Introducing $\tilde{w}\in\Sigma_K$ defined by
\[\forall a = 1 \dots K, \ \ \tilde{w}_a = \frac{1/\Delta_a^2}{\sum_{i=1}^K 1/\Delta_i^2}\;,\]
it holds that 
\[T^*(\bm\mu)^{-1} \geq \inf_{a \neq 1} \frac{\tilde{w}_1\tilde{w}_a}{\tilde{w}_1 + \tilde{w}_a}\frac{\Delta_a^2}{2\sigma^2} = \frac{1}{\sum_{i=1}^K\frac{2\sigma^2}{\Delta_i^2}}\inf_{a \neq 1}\frac{1}{1 + \frac{\Delta_1^2}{\Delta_a^2}}.\]
The infimum is obtained for $a=2$, and using that $\Delta_2=\Delta_1$ leads to the upper bound
\[T^*(\bm\mu) \leq 2 \left(\sum_{i=1}^K \frac{2\sigma^2}{\Delta_i^2}\right)\;.\]
 The following lower bound was obtained by \cite{JMLR15} for every PAC strategy: 
\[\liminf_{\delta \rightarrow 0}\frac{\bE[\tau_\delta]}{\log(1/\delta)} \geq \sum_{i=1}^K \frac{2\sigma^2}{\Delta_i^2}\;.\]
Combining this inequality with the upper bound on $\liminf_{\delta \rightarrow 0}{\bE[\tau_\delta]}/{\log(1/\delta)}$ obtained in Theorem~\ref{thm:AsymptoticSC} for the $\delta$-PAC Track-and-Stop algorithm shows that 
\[T^*(\bm\mu) \geq \sum_{i=1}^K \frac{2\sigma^2}{\Delta_i^2}\;,\]
which concludes the proof.

\section{Tracking results}\label{proof:Tracking}

Lemma~\ref{lem:TrackingC} and \ref{lem:TrackingP} both follow from deterministic results that we can give for procedures tracking any cumulated sums of proportions (Lemma~\ref{lem:pureTracking}) or any changing sequence of proportions that concentrates (Lemma~\ref{lem:trackingCsaba}). We state and prove theses two results in this section, and also explain how they lead to Lemma~\ref{lem:TrackingC} and \ref{lem:TrackingP}. We then provide a proof of Proposition~\ref{prop:ConvergenceFraction}. 

\subsection{Tracking a cumulated sum of proportions}
\begin{lemma}\label{lem:pureTracking}
	Let $K$ be a positive integer, let $\Sigma_K$ be the simplex of dimension $K-1$, and for every $i\in\{1,\dots,K\}$, let $\delta_i$ be the vertex of $\Sigma_K$ with a $1$ on coordinate $i$.
	For a positive integer $n$, let $p(1), p(2),\dots, p(n)\in \Sigma_K$ and for every $k\leq n$ let $P(k) = p(1)+\dots+p(k)$.	
	Define $N(0)=0$, for every $k\in\{0,\dots,n-1\}$
	 \[I_{k+1} \in \argmax{1\leq i\leq K} \ \left[ P_{i}(k+1)-N_{i}(k)\right]\] and 
	$N({k+1}) = N({k}) + \delta_{I_{k+1}}$. Then
	\[\max_{1\leq i\leq K} \big|N_i(n)-P_i(n)\big| \leq K-1\;. \]
\end{lemma}
To obtain Lemma~\ref{lem:TrackingC}, we start by applying Lemma~\ref{lem:pureTracking} with $p(k) = w^{\epsilon_{k-1}}(\hat{\bm \mu}(k-1))$, so that $P({k+1}) = \sum_{s=0}^{k} w^{\epsilon_s}(\hat{\bm \mu}(s))$. One obtains
\begin{equation}\max_{1 \leq a \leq K} \left|N_a(t) - \sum_{s=0}^{t-1} w^{\epsilon_s}_a(\hat{\bm\mu}(s)) \right| \leq K-1\;.\label{mainIneq}\end{equation}
Moreover, by definition of $w^{\epsilon}(s)$,  
\[\max_{1 \leq a \leq K} \left|\sum_{s=0}^{t-1} w^{\epsilon_s}_a(\hat{\bm\mu}(s)) - \sum_{s=0}^{t-1} w^{*}_a(\hat{\bm\mu}(s))  \right| \leq \sum_{s=0}^{t-1} K \epsilon_s\;.\]
Now, with the choice $\epsilon_t = (K^2+t)^{-1/2}/2$, one has 
\[\sqrt{t+K^2} - K = \int_{0}^{t} \frac{ds}{2\sqrt{K^2+s}} \leq \sum_{s=0}^{t-1} \epsilon_s \leq \int_{-1}^{t-1} \frac{ds}{2\sqrt{K^2+s}} = \sqrt{t + K^2 -1} - \sqrt{K^2-1}\;,\]
which yields 
\[\max_{1 \leq a \leq K} \left|N_a(t) - \sum_{s=0}^{t-1} w^{*}_a(\hat{\bm\mu}(s)) \right| \leq K-1 + K\left(\sqrt{t + K^2 -1} - \sqrt{K^2-1}\right) \leq K(1 + \sqrt{t})\;.\]
From \eqref{mainIneq}, it also follows that 
\[N_a(t) \geq \sum_{s=0}^{t-1} \epsilon_s - (K-1) \geq  \sqrt{t+K^2} -K - (K-1)  \geq \sqrt{t+K^2} - 2K,\]
which concludes the proof of Lemma~\ref{lem:TrackingC}. 

\paragraph{Proof of Lemma~\ref{lem:pureTracking}.} First, we prove by induction on $k$ that  
\[\max_{1 \leq i \leq K} N_i(k) - P_i(k) \leq 1. \]
The statement is obviously true for $k=0$. Assume that it holds for some $k \geq 0$. For $i \neq I_{k+1}$ one has 
$N_{i}(k+1) - P_{i}(k+1) = N_i(k) - P_{i}(k) - p_i(k) \leq 1 - p_i(k) \leq 1$, whereas
$N_{I_{k+1}}(k+1) - P_{I_{k+1}}(k+1) = 1 + (N_{I_{k+1}}(k) - P_{I_{k+1}}(k+1) \leq 1$, using that $\sum_{i}(P_{i}(k+1) - N_i(k))=0$, hence the largest term in this sum (which is for $i=I_{k+1}$ by definition) is positive. 

It follows that, for all $k$,  
	\begin{align*}\max_{1\leq i\leq K} \big|N_i(k)-P_i(k)\big| & = \max\big\{ \max_{1\leq i\leq K} P_i(k)-N_i(k) , \max_{1\leq i\leq K} N_i(k)-P_i(k) \big\}\\
	& \leq \max\{ \sum_{1\leq i\leq K}\big(P_i(k)-N_i(k) \big)_+  , 1 \}\;.
	\end{align*}
Introducing, for every $k\in\{1,\dots,n\}$
\[r_k = \sum_{i=1}^K\big(P_i(k)-N_i(k) \big)_+,\]
Lemma~\ref{lem:pureTracking} follow from the following bound on $r_k$, that we prove by induction on $k$:
\[r_k \leq K-1.\]
Observe that $r_0=0$. For every $k\geq 0$, one can write
\[r_{k+1} = r_k + \sum_{i = 1}^K  p_{i}(k+1)\ind_{(P_{i}(k+1) \geq N_{i}(k+1))} - \ind_{(P_{I_{k+1}}(k+1) \geq N_{I_{k+1}}(k+1))}.\]
We distinguish two cases. If $(P_{I_{k+1}}(k+1) \geq N_{I_{k+1}}(k+1))$ one has \[r_{k+1} \leq r_k +  \sum_{i=1}^K  p_{i}(k+1) - 1 = r_k \leq K-1.\] If $(P_{I_{k+1}}(k+1) < N_{I_{k+1}}(k+1))$, which implies $P_{I_{k+1}}(k+1) - N_{I_{k+1}}(k) \leq 1$, one has
\begin{eqnarray*}
r_{k+1} &=& \sum_{\substack{i =1  \\i \neq I_{k+1} }}^K (P_{i}(k+1) - N_{i}(k+1))_+ \leq \sum_{\substack{i =1  \\i \neq I_{k+1} }}^K \max_j \left(P_{j}(k+1) - N_{j}(k)\right)\\
&{=} & \sum_{\substack{i =1  \\i \neq I_{k+1} }}^K \left(P_{I_{k+1}}(k+1) - N_{I_{k+1}}(k)\right) \leq (K-1)\;,
\end{eqnarray*}
which concludes the proof.

\begin{remark} This result is probably overly pessimistic, but one cannot hope for an upper bound independent of $K$: for $p_i(k) = 1\{i\geq k\}/(K-k+1)$, by choosing at each step $K$ the smallest index in the $\argmax{}$ defining $I_k$, one gets $N_{K}(K-1)=0$ and $P_{K}(K-1)=1/K+1/(K-1)+...+1/2\sim \log(K)$.
\end{remark}

\subsection{Tracking a changing sequence that concentrates}

\begin{lemma}\label{lem:trackingCsaba}
	Let $K$ be a positive integer and let $\Sigma_K$ be the simplex of dimension $K-1$.
	Let $g:\N\to\R$ be a non-decreasing function such that $g(0)=0$, $g(n)/n \rightarrow 0$ when $n$ tends to infinity 
	and for all $k\geq 1$ and $\forall m\geq 1$,
    \[\inf\{k\in \N : g(k)\geq m\} > \inf\{k\in \N: g(k)\geq m-1\} +K.\]	
    	Let $\hat{\lambda}(k)$ be a sequence of elements in $\Sigma_K$ such that there exists $\lambda^* \in \Sigma_K$, there exists $\epsilon>0$ and an integer $n_0(\epsilon)$ such that  
    	\[\forall n \geq n_0 , \ \sup_{1 \leq i \leq K} | \hat{\lambda}_i(k) - \lambda^*_i| \leq \epsilon\;.\]

	Define $N(0)=0$, and for every $k\in\{0,\dots,n-1\}$, $U_k = \{ i : N_i(k) < g(k)\}$ and 
	\[I_{k+1} \in \left\{\begin{array}{l}
	                    \underset{i \in U_k}{\emph{argmin}} \ N_{i}(k) \ \ \ \  \text{if} \ U_k \neq \emptyset, \\
	                    \underset{i \in \{1,\dots,K\}}{\emph{argmax}} \ \left[k\hat{\lambda}_i(k) - N_i(k)\right] \ \ \text{else},
	                     \end{array}
        \right.
        \]
	and  for all $i$ $N_{i}(k+1) = N_{i}(k) + \ind_{(I_{k+1}=i)}$.
       Then for all $i \in \{1,\dots,K\}$, $N_i(n) > g(n) -1$ and there exists $n_1 \geq n_0$ (that depends on $\epsilon$) such that for all $n \geq n_1$, 
    \[\max_{1\leq i\leq K} \left|\frac{N_n(i)}{n}-\lambda^*_i\right| \leq 3(K-1)\epsilon\;.\]
\end{lemma}

First it is easy to check that $g(k)=(\sqrt{k} - K/2)_+$ satisfies the assumptions of Lemma~\ref{lem:trackingCsaba}. Lemma~\ref{lem:TrackingP} then  follows by choosing $\hat{\lambda}(k) = w^*(\hat{\mu}(k))$ and $\lambda^* = w^*(\bm \mu)$. The constant $n_1$ in Lemma~\ref{lem:trackingCsaba} depends on $\epsilon$, hence the notation $t_\epsilon$. 

The proof of this Lemma~\ref{lem:trackingCsaba} is inspired by the proof of Lemma 3 in \cite{CsabaTracking}, although a different tracking procedure is analyzed. 

\paragraph{Proof of Lemma~\ref{lem:trackingCsaba}.} First, we justify that $N_i(n) > g(n)-1$. For this purpose, we introduce for all $m\in \N$ the integer $n_m=\inf\{k \in \N : g(k)\geq m\}$. We also let $\cI_m = \{ n_m, \dots, n_{m+1} -1 \}$. From our assumption on $g$, it follows that $|\cI_m| > K$ and by definition, for all $k \in \cI_m$, $m \leq g(k) < m +1$. 

We prove by induction that the following statement holds for all $m$: 
\begin{equation}\forall k \in \cI_m, \forall i,  N_i(k) \geq m.\ \ \text{Moreover} \ \text{for} \ k \geq n_{m} + K,  \ U_k = \emptyset \ \ \text{and} \ N_i(k) \geq m+1\;.\label{HR1} \end{equation}
First, for all $k \in \cI_0$, one has $U_k = \{ i : N_i(k)=0 \}$. Therefore $(I_1,\dots,I_K)$ is a permutation of $(1,\dots,K)$, thus for  $k \geq K = n_0 + K$, $N_i(k) \geq 1$, and $U_k = \emptyset$, and the statement holds for $m=0$. Now let $m\geq 0$ such that the statement is true. From the inductive hypothesis, one has 
 \[\forall k \in \cI_{m+1}, \forall i,   N_{i}(k) \geq N_i(n_{m+1}-1) \geq m + 1\;.\] 
Besides, as $g(k) < m+2$ for $k \in \cI_{m+1}$, one has $U_k = \{ i : N_i(k) = m +1\}$ and $I_k$ is chosen among this set while it is non empty. For $k \geq n_{m+1} +K$, it holds that $U_k = \emptyset$ and $N_i(k) \geq m +2$ for all $i$. Thus the statement holds for $m+1$.  

From the fact that \eqref{HR1} holds for all $m$, using that for $k \in \cI_m$, $m > g(k) -1$, it follows that for all $k$, for all $i$,  $N_i(k) > g(k) -1 $.

Now for all $i \in \{1,\dots,K\}$, we introduce $E_{i,n} = N_{i}(n) - n \lambda^*_i$. Using that 
\begin{equation}\sum_{i=1}^K E_{i,n} = 0,\label{Nul}\end{equation}
leads to
\[\sup_i |E_{i,n}| \leq (K-1) \sup_i \ E_{i,n}.\]
Indeed, for every $i$, one has $E_{i,n} \leq \sup_{k} E_{k,n}$ and  
\[E_{i,n} = - \sum_{j \neq i} E_{j,n} \geq - \sum_{j \neq i} \sup_{k} E_{k,n} = - (K-1)\sup_{k}E_{k,n}\;.\] 

To conclude the proof, we give an upper bound on $\sup_i E_{i,n}$, for $n$ large enough. Let $n_0' \geq n_0$ such that 
\[\forall n \geq n_0', \ \ \ g(n) \leq 2n \epsilon \ \ \text{and} \ \ 1/n \leq \epsilon\;.\]
We first show that for $n\geq n_0'$,  
\begin{equation}
(I_{n+1} = i) \subseteq \left(E_{i,n} \leq 2n \epsilon\right)\label{IntermediateEvent}
\end{equation}
To prove this, we write 
\[(I_{n+1} = i) \subseteq \left(N_i(n) \leq g(n) \ \ \text{or} \ \ i = \underset{1 \leq j \leq K}{\emph{argmin}} \ \left(N_j(n) - n \hat{\lambda}_j(n)\right)\right)\;.
\]
Now if $N_i(n) \leq g(n)$, one has $E_{i,n} \leq g(n) - n\lambda^*_i \leq g(n) \leq 2n\epsilon$,  by definition of $n_0'$.  In the second case, one has 
\begin{eqnarray*}
 N_i(n) - n \hat{\lambda}_i(n) &= &\min_j \left(N_j(n) - n \hat{\lambda}_j(n)\right) \\
 E_{i,n} + n(\lambda^*_i - \hat{\lambda}_i(n)) & = & \min_j  \left(E_{j,n} + n(\lambda^*_j - \hat{\lambda}_j(n)\right) 
\end{eqnarray*}
Using the closeness of each $\hat{\lambda}_j(n)$ to the corresponding $\lambda_j^*$, as $n\geq n_0$, yields
\[E_{i,n} + n(\lambda^*_i - \hat{\lambda}_n(i)) \leq  \min_j  \left(E_{j,n} + n\epsilon\right) \leq n\epsilon,\]
where we use that $\min_j E_{j,n} \leq 0$ by \eqref{Nul}. Using that $|\lambda^*_i - \hat{\lambda}_i(n)| \leq \epsilon$ as well, one obtains 
\[E_{i,n} \leq 2n\epsilon\;.\]
This proves \eqref{IntermediateEvent}. 

$E_{i,n}$ satisfies $E_{i,n+1} = E_{i,n} + \ind_{(I_{n+1}=i)} - \lambda_i^*$, therefore, if $n\geq n_0'$, 
\begin{eqnarray*}
 E_{i,n+1} &\leq& E_{i,n} + \ind_{(E_{i,n} \leq 2n\epsilon)} - \lambda_i^*\;. 
\end{eqnarray*}
We now prove by induction that for every $n \geq n_0'$, one has 
\[E_{i,n} \leq \max(E_{i,n_0'}, 2n\epsilon +1 ).\]
For $n=n_0'$, this statement clearly holds. Let $n\geq n_0'$ such that the statement holds.
If $E_{i,n} \leq 2n\epsilon$, one has 
\begin{eqnarray*}E_{i,n+1} &\leq& 2n\epsilon + 1 - \lambda_i^* \leq 2n\epsilon +1 \leq \max(E_{i,n_0'},2(n)\epsilon+1) \\
&\leq& \max(E_{i,n_0'},2(n+1)\epsilon +1).\end{eqnarray*}
If $E_{i,n} > 2n\epsilon$, the indicator is zero and  
\[E_{i,n+1} \leq \max(E_{i,n_0'},2n\epsilon +1) - \lambda_i^* \leq \max(E_{i,n_0'},2(n+1)\epsilon+1)\;,\]
which concludes the induction. 

For all $n\geq n_0'$, using that $E_{i,n_0'} \leq n_0'$ and $1/n \leq \epsilon$, it follows that 
\[\sup_i \left|\frac{E_{i,n}}{n}\right| \leq  (K-1) \max\left(2\epsilon + \frac{1}{n}, \frac{n_0'}{n}\right) \leq (K-1) \max\left(3\epsilon, \frac{n_0'}{n}\right).\]
Hence there exists $n_1 \geq n_0'$ such that, for all $n\geq n_1$, 
\[\sup_i \left|\frac{E_{i,n}}{n}\right| \leq  3(K-1) \epsilon\;,\]
which concludes the proof.

\subsection{Proof of Proposition~\ref{prop:ConvergenceFraction}}\label{proof:ConvergenceFraction}

Because of the forced-exploration step, both tracking strategies satisfy $\forall a, N_a(t) \rightarrow \infty$. Thus, from the law of large number, the event 
\[\cE = \left(\hat{\bm\mu}(t) \underset{t \rightarrow \infty}{\rightarrow} \bm \mu\right)\]
is of probability one. 
For C-Tracking, it follows from Lemma~\ref{lem:TrackingC} that
\[\left|\frac{N_a(t)}{t} - \frac{1}{t}\sum_{s=0}^{t-1} w_a^*(\bm \mu(s))\right| \leq \frac{K(\sqrt{t} +1)}{t}\;.\]
By continuity of $w^*$ in $\bm \mu$, for all $a$, $w_a^*(\hat{\bm\mu}(t)) \rightarrow w^*_a(\bm\mu).$ Using moreover the Cesaro lemma, one obtains that, on $\cE$, $N_a(t)/t \rightarrow w_a^*(\bm\mu)$.
For D-Tracking, we first use that for $\omega \in \cE$, there exists $t_0(\epsilon)$ such that 
\[\sup_{t \geq t_0(\epsilon)} \max_ a |w_a^*(\hat{\bm\mu}(t)) - w_a^*(\bm \mu)| \leq \frac{\epsilon}{3(K-1)}\;,\]
by continuity of the function $\bm \lambda \mapsto w^*(\bm \lambda)$ in $\bm \mu$. Hence, using Lemma~\ref{lem:TrackingP}, there exists $t_\epsilon \geq t_0$ such that for all $t \geq t_\epsilon$, 
\[\max_a\left| \frac{N_a(t)}{t} - w_a^*(\bm \mu) \right| \leq \epsilon\;.\]
Hence, for this $\omega \in \cE$,  $N_a(t) / t (\omega) \rightarrow w_a^*(\bm \mu)$ for all $a$.

\section{PAC guarantees}\label{proof:deltaPAC}

\subsection{Proof of Proposition~\ref{prop:deltaPAC}.}\label{sec:proof:deltaPAC}
Recall that 
\begin{align*}
\bP_{\bm \mu}(\tau_\delta<\infty,\hat{a}_{\tau_\delta}\neq a^*) & \leq \bP_{\bm \mu}\left(\exists a \in \Arms\setminus a^*, \exists t \in \N : Z_{a,a^*}(t) > \beta(t,\delta)\right) \\& \leq \sum_{a\in \Arms\setminus A^*} \bP_{\bm \mu}(T_{a,a^*} < \infty)\;,
\end{align*}
where $T_{a,b} := \inf\{ t\in \N : Z_{a,b}(t) > \beta(t,\delta)\}$. To conclude the proof, we now show that for $\beta(t,\delta)$ as in Proposition~\ref{prop:deltaPAC}, for any $a,b$ such that $\mu_a<\mu_b$, one has 
\[\bP_{\bm \mu}(T_{a,b} < \infty) \leq \frac{\delta}{K-1}.\]

Let $a$,$b$ be such that $\mu_a < \mu_b$. One introduces $f_{\bm \mu}(\underline{x}_t,\underline{a}_t)$ the likelihood of observing the sequence of rewards $\underline{x}_t = (x_1,\dots,x_t)$ and sequence of actions $\underline{a}_t=(a_1,\dots,a_t)$. One has 
\[f_{\bm \mu}(\underline{x}_t,\underline{a}_t) = \prod_{i=1}^{K}p_{\mu_i}\left(\underline{x}_t^i\right)\times \left[p(a_1) \prod_{s=2}^t p(a_s | \underline{x}_{s-1},\underline{a}_{s-1})\right],\]
where $\underline{x}_t^i$ is a vector that gathers the sequence of successive rewards from arm $i$ up to time $t$ (that is a function of both $\underline{x}_t$ and $\underline{a}_t$). $f_{\bm \mu}(\underline{x}_t,\underline{a}_t)$ is a probability density on $\cX^t\times \Arms^t$. For any density $h$ supported on $\cI$, the (partially) integrated likelihood
\[I_h(\underline{x}_t,\underline{a}_t) = \prod_{i\in\Arms \setminus \{a,b\}}p_{\mu_i}\left(\underline{x}_t^i\right) \left(\int_\R p_{u}\left(\underline{x}_t^a\right)h(u)du\right)\left(\int_\R p_{u}\left(\underline{x}_t^b\right)h(u)du\right)\times \left[p(a_1) \prod_{s=2}^t p(a_s | \underline{x}_{s-1},\underline{a}_{s-1})\right]\]
is also a probability distribution.

On the event $(T_{a,b}=t)$, $Z_{a,b}(t)$ exceeds for the first time the threshold $\beta(t,\delta)$, which implies in particular (from the definition of $Z_{a,b}(t)$) that  
\[1 \leq e^{-\beta(t,\delta)} \frac{\max_{{\mu}_a' \geq {\mu}'_b} p_{\mu_a'}(\underline{X}^a_t)p_{\mu_b'}(\underline{X}^b_t)}{\max_{{\mu}_a' \leq {\mu}'_b} p_{\mu_a'}(\underline{X}^a_t)p_{\mu_b'}(\underline{X}^b_t)}\;.\]
We use this fact in the first inequality below, whereas the second inequality is based on the fact that $\bm \mu$ satisfies $\mu_a < \mu_b$:
 \begin{align*}
 \bP_{\bm \mu}(T_{a,b}< \infty) & = \sum_{t=1}^\infty \bP_{\bm \mu}(T_{a,b}=t) = \sum_{t=1}^\infty \bE_{\bm \mu}\left[\ind_{(T_{a,b}=t)}\right] \\
 & \leq \sum_{t=1}^\infty e^{-\beta(t,\delta)}\bE_{\bm \mu}\left[\ind_{(T_{a,b}=t)} \frac{\max_{{\mu}_a' \geq {\mu}'_b} p_{\mu_a'}(\underline{X}^a_t)p_{\mu_b'}(\underline{X}^b_t)}{\max_{ {\mu}_a' \leq {\mu}'_b} p_{\mu_a'}(\underline{X}^a_t)p_{\mu_b'}(\underline{X}^b_t)}\right] \\
 & \leq \sum_{t=1}^\infty e^{-\beta(t,\delta)}\underbrace{\bE_{\bm \mu}\left[\ind_{(T_{a,b}=t)} \frac{\max_{{\mu}_a' \geq {\mu}'_b} p_{\mu_a'}(\underline{X}^a_t)p_{\mu_b'}(\underline{X}^b_t)}{p_{\mu_a}(\underline{X}^a_t)p_{\mu_b}(\underline{X}^b_t)}\right]}_{(*)}\;. 
\end{align*}
We now explicit the expectation $(*)$ for Bernoulli distributions. 
\begin{align*}
(*) & = \sum_{\underline{a}_t \in \Arms^t}\sum_{\underline{x}_t \in \{0,1\}^t} \ind_{(T_{a,b}=t)}(\underline{x}_t,\underline{a}_t) \max_{{\mu}_a' \geq {\mu}'_b}p_{\mu_a'}(\underline{x}^a_t)p_{\mu_b'}(\underline{x}^b_t) \prod_{i \in \Arms \setminus \{a,b\}} p_{\mu_i}(\underline{x}_t^i) \left[p(a_1) \prod_{s=2}^t p(a_s | \underline{x}_{s-1},\underline{a}_{s-1})\right]
\end{align*}
Introducing, for a vector $x$, 
\[\kt(x) = \int_0^1 \frac{1}{\sqrt{\pi u(1-u)}} p_u(x) \d u\;,\]
we have that 
\[ \max_{{\mu}_a' \geq {\mu}'_b}p_{\mu_a'}(\underline{x}^a_t)p_{\mu_b'}(\underline{x}^b_t) \leq 2 t \times  \kt(\underline{x}^a_t)\kt(\underline{x}^b_t),\]
which follows from Lemma~\ref{lem:ktMain} stated in the main text and the fact that for $n_a,n_b$ are such that $n_a + n_b \leq t$, $4\sqrt{n_a}\sqrt{n_b}\leq 2(n_a + n_b) \leq 2t$.
Using this inequality to upper bound $(*)$, one recognize the integrated likelihood associated to the density $h(u)=\frac{1}{\sqrt{\pi u (1-u)}}\ind_{[0,1]}(u)$:
 \begin{align*}
(*)  & \leq 2t \sum_{\underline{a}_t \in \Arms^t}\sum_{\underline{x}_t \in \{0,1\}^t} \ind_{(T_{a,b}=t)}(\underline{x}_t,\underline{a}_t) I_{h}(\underline{x}_t,\underline{a}_t) =  2t \tilde{\bP}(T_{a,b}=t)\;, 
\end{align*}
where $\tilde{\bP}$ is an alternative probabilistic model, under which $\mu_a$ and $\mu_b$ are drawn from a $\mathrm{Beta}(1/2,1/2)$ (prior) distribution at the beginning of the bandit game. Finally, using the explicit expression of $\beta(t,\delta)$,
\begin{eqnarray*}
 \bP_{\bm \mu}(T_{a,b}< \infty) & \leq & \sum_{t=1}^\infty 2t e^{-\beta(t,\delta)} \tilde{\bP}(T_{a,b}=t) \leq \frac{\delta}{K-1} \sum_{t=1}^\infty \tilde{\bP}(T_{a,b}=t) \\ &= &\frac{\delta}{K-1} \tilde{\bP}(T_{a,b}< \infty) \leq \frac{\delta}{K-1}\;,
\end{eqnarray*}
which concludes the proof.

\subsection{Proof of Proposition~\ref{prop:PACGene}.}\label{proof:PACGene} The proof relies on the fact that $Z_{a,b}(t)$ can be expressed using function $I_{\alpha}$ introduced in Definition~\ref{def:Divergences}. An interesting property of this function, that we use below, is the following. It can be checked that if $x > y$, 
\[I_\alpha(x,y) = \inf_{x'<y'} \left[\alpha d(x,x') + (1-\alpha)d(y,y')\right].\]
For every $a,b$ that are such that $\mu_a<\mu_b$ and $\hat{\mu}_a(t)>\hat{\mu}_b(t)$, one has the following inequality: 
\begin{eqnarray*}
Z_{a,b}(t) &=& (N_a(t) + N_b(t)) I_{\frac{N_a(t)}{N_a(t) + N_b(t)}}\left(\hat{\mu}_a(t),\hat{\mu}_b(t)\right) \\
 &= &\inf_{\mu_a'<\mu_b'}{N_a(t)}d(\hat{\mu}_{a}(t),\mu_a') + {N_b(t)}d(\hat{\mu}_{b}(t),\mu_b')  \\
&\leq & {N_a(t)d\left(\hat{\mu}_{a}(t),\mu_a\right) + N_b(t)d(\hat{\mu}_b(t),\mu_b)}\;.
\end{eqnarray*}
One has
\begin{align*}
 & \bP_{\bm\mu}(\tau_\delta<\infty,\hat{a}_{\tau_\delta}\neq a^*)   \leq   \bP_{\bm\mu}\left(\exists a \in \Arms\setminus a^*, \exists t \in \N : \hat{\mu}_a(t)> \hat{\mu}_{a^*}(t), Z_{a,a^*}(t) > \beta(t,\delta)\right) \\
& \hspace{2cm}\leq \bP_{\bm\mu}\big(\exists t\in \N : \exists a\in \Arms\setminus a^* :  N_a(t)d\left(\hat{\mu}_{a}(t),\mu_a\right) + N_{a^*}(t)d(\hat{\mu}_{a^*}(t),\mu_{a^*}) \geq \beta(t,\delta)\big)\\
& \hspace{2cm}\leq \bP_{\bm\mu}\left(\exists t\in \N : \sum_{a=1}^K N_a(t) d(\hat{\mu}_a(t),\mu_a) \geq \beta(t,\delta)\right) \\
& \hspace{2cm}\leq \sum_{t=1}^\infty e^{K+1}\left(\frac{\beta(t,\delta)^2\log(t)}{K}\right)^Ke^{-\beta(t,\delta)}\;.
\end{align*}
The last inequality follows from a union bound and Theorem 2 of \cite{Combes14Lip}, originally stated for Bernoulli distributions but whose generalization to one-parameter exponential families is straightforward. Hence, with an exploration rate of the form $\beta(t,\delta)=\log(C t^\alpha / \delta)$, for $\alpha >1$, choosing $C$ satisfying
\[\sum_{t=1}^\infty \frac{e^{K+1}}{K^K} \frac{(\log^2(Ct^\alpha) \log t)^K}{t^\alpha} \leq C\]
yields a probability of error upper bounded by $\delta$.

\section{Expected sample complexity analysis}\label{sec:ProofSC}

The proof of Theorem~\ref{thm:AsymptoticSC} relies on two ingredients: a concentration result for the empirical mean $\hat{\bm\mu}(t)$, that follows from the forced exploration (Lemma~\ref{lem:concSimple}) and the tracking lemma associated to the sampling strategy used (Lemma~\ref{lem:TrackingC}, Lemma~\ref{lem:TrackingP}).  
We start with a small technical lemma that can be checked directly, or that can be seen as a by-product of well-known bounds on the Lambert $W$ function.
\begin{lemma}\label{lem:technical} For every $\alpha \in [1,e/2]$, for any two constants $c_1,c_2>0$, 
\[x = \frac{\alpha}{c_1}\left[\log\left(\frac{c_2 e}{c_1^\alpha}\right) + \log\log\left(\frac{c_2}{c_1^\alpha}\right)\right]\]
is such that $c_1  x \geq \log(c_2 x^\alpha)$. 
\end{lemma}
\subsection{Proof of Theorem~\ref{thm:AsymptoticSC}}

To ease the notation, we assume that the bandit model $\bm \mu$ is such that $\mu_1 > \mu_2 \geq \dots \geq \mu_K$.
Let $\epsilon>0$. From the continuity of $w^*$ in $\bm \mu$, there exists $\xi=\xi(\epsilon) \leq (\mu_1-\mu_2)/4$ such that 
\[\cI_\epsilon := [\mu_1 - \xi,\mu_1 + \xi] \times \dots \times [\mu_K - \xi, \mu_K + \xi]\]
is such that for all $\bm\mu' \in \cI_\epsilon$, \[\max_a |w_a^*(\bm\mu') - w^*_a(\bm\mu)| \leq \epsilon.\] In particular, whenever $\hat{\bm \mu}(t) \in \cI_\epsilon$, the empirical best arm is $\hat{a}_t =1$.

Let $T \in \N$ and define $h(T):=T^{1/4}$ and the event 
\[\cE_T(\epsilon)= \bigcap_{t = h(T)}^{T}\left(\hat{\bm \mu}(t) \in \cI_\epsilon \right).\]
The following lemma is a consequence of the ``forced exploration'' performed by the algorithm, that ensures that each arm is drawn at least of order $\sqrt{t}$ times at round $t$. 

\begin{lemma}\label{lem:concSimple} There exist two constants $B,C$ (that depend on $\bm \mu$ and $\epsilon$) such that \[\bP_{\bm \mu}(\cE_T^c) \leq B T \exp(-C T^{1/8}).\] 
\end{lemma}

Then, exploiting the corresponding tracking Lemma, one can prove the following 

\begin{lemma}\label{lem:CcqceTracking} There exists a constant $T_\epsilon$ such that for $T \geq T_\epsilon$, it holds that on $\cE_T$, for either C-Tracking or D-Tracking, 
 \[\forall t \geq \sqrt{T}, \ \max_a \left|\frac{N_a(t)}{t} - w_a^*(\mu)\right| \leq 3(K-1)\epsilon\]
\end{lemma}

\paragraph{Proof.} This statement is obvious for D-Tracking, just by definition of $\cI_\epsilon$ and by Lemma~\ref{lem:TrackingP}. For C-Tracking, for any $t \geq \sqrt{T}=h(T)^2$, using Lemma~\ref{lem:TrackingC}, one can write, for all $a$,
\begin{eqnarray*}
 \left|\frac{N_a(t)}{t} - w_a^*(\bm \mu)\right| & \leq & \left|\frac{N_a(t)}{t} - \frac{1}{t} \sum_{s=0}^{t-1} w_a^*(\hat{\bm \mu}(s))\right| + \left|\frac{1}{t} \sum_{s=0}^{t-1} w_a^*(\hat{\bm \mu}(s)) - w_a^*(\bm \mu)\right| \\
 & \leq & \frac{K(\sqrt{t} +1)}{{t}} + \frac{h(T)}{t} + \frac{1}{t}\sum_{t=h(T)}^{t-1} |w_a^*(\hat{\bm \mu}(s)) - w_a^*(\bm \mu)| \\
 &\leq &\frac{2K}{T^{1/4}} + \frac{1}{h(T)} + \epsilon = \frac{2K + 1 }{T^{1/4}} + \epsilon \leq 3\epsilon\;,
\end{eqnarray*}
whenever $T \geq ((2K+1)/2\epsilon)^4$. 
\qed

On the event $\cE_T$, it holds for $t \geq h(T)$ that $\hat{a}_t=1$ and the Chernoff stopping statistic rewrites 
\begin{eqnarray*} 
 \max_a \min_{b \neq a} \ Z_{a,b}(t) & = & \min_{a \neq 1} \ Z_{1,a}(t) =  \min_{a \neq 1} N_1(t)d(\hat{\mu}_1(t) , \hat{\mu}_{1,a}(t)) + N_a(t) d(\hat{\mu}_a(t),\mu_{1,a}(t))\\
 & = & t \left[\min_{a \neq 1} \left(\frac{N_1(t)}{t} + \frac{N_a(t)}{t}\right)I_{\frac{N_1(t)/t}{N_1(t)/t + N_a(t)/t}}\left(\hat{\mu}_1(t),\hat{\mu}_a(t)\right)  \right] \\
 & = & t  g \left(\hat{\bm \mu}(t),\left(\frac{N_a(t)}{t}\right)_{a=1}^K\right)\;, 
\end{eqnarray*}
where we introduce the function 
\[g(\bm \mu',\bm w') = \min_{a \neq 1}(w_1' + w'_a) I_{\frac{w'_1}{w'_1 + w'_a}}(\mu_1',\mu_a').\]
Using Lemma~\ref{lem:CcqceTracking}, for $T\geq T_\epsilon$, introducing
\[C^*_\epsilon(\bm \mu) = \inf_{\substack{\bm \mu' : ||\bm\mu' - \bm \mu|| \leq \alpha(\epsilon) \\ \bm w' : ||\bm w' -w^*(\bm \mu)||\leq 2(K-1)\epsilon}} g(\bm \mu',\bm w')\;,\]
on the event $\cE_T$ it holds that for every $t \geq \sqrt{T}$, 
\[\left(\max_a \min_{b \neq a} \ Z_{a,b}(t) \geq t C^*_\epsilon(\bm \mu)\right)\;.\]

Let $T \geq T_\epsilon$. On $\cE_T$, 
\begin{eqnarray*}
 \min(\tau_\delta,T) & \leq & \sqrt{T} + \sum_{t=\sqrt{T}}^T \ind_{\left(\tau_\delta > t\right)} \leq \sqrt{T} + \sum_{t=\sqrt{T}}^T \ind_{\left(\max_a \min_{b \neq a} \ Z_{a,b}(t) \leq \beta(t,\delta)\right)} \\
 & \leq & \sqrt{T} + \sum_{t=\sqrt{T}}^T \ind_{\left(t C_\epsilon^*(\bm\mu) \leq \beta(T,\delta)\right)} \leq \sqrt{T} + \frac{\beta(T,\delta)}{C_\epsilon^*(\bm \mu)}\;.
\end{eqnarray*}
Introducing 
\[T_0(\delta) = \inf \left\{ T \in \N : \sqrt{T} + \frac{\beta(T,\delta)}{C_\epsilon^*(\bm \mu)} \leq T \right\},\]
for every $T \geq \max (T_0(\delta), T_\epsilon)$, one has $\cE_T \subseteq (\tau_\delta \leq T)$, therefore 
\[\bP_{\bm \mu}\left(\tau_\delta > T\right) \leq \bP(\cE_T^c) \leq BT \exp(-C T^{1/8})\]
and
\[\bE_{\bm \mu}[\tau_\delta] \leq T_0(\delta) + T_\epsilon + \sum_{T=1}^\infty BT \exp(-C T^{1/8})\;.\]
We now provide an upper bound on $T_0(\delta)$. Letting $\eta >0$ and introducing the constant 
\[C(\eta) = \inf \{ T \in \N : T - \sqrt{T} \geq T/(1+\eta)\}\]
one has 
\begin{eqnarray*}
 T_0(\delta) & \leq & C(\eta) + \inf \left\{T \in \N : \frac{1}{C_\epsilon^*(\bm\mu)} \log\left(\frac{r(T)}{\delta}\right) \leq \frac{T}{1+\eta}\right\} \\
 & \leq & C(\eta) + \inf \left\{T \in \N : \frac{C_\epsilon^*(\bm\mu)}{1+\eta} T \geq \log\left(\frac{Dt^{1+\alpha}}{\delta}\right) \right\},
\end{eqnarray*}
where the constant $D$ is such that $r(T) \leq D T^\alpha$. Using again the technical Lemma~\ref{lem:technical}, one obtains, for $\alpha \in [1,e/2]$, 
\[T_0(\delta) \leq C(\eta) + \frac{\alpha(1+\eta)}{C_\epsilon(\bm \mu)}\left[\log \left(\frac{De(1+\eta)^\alpha}{\delta (C_\epsilon^*(\bm\mu))^\alpha}\right) +  \log\log \left(\frac{D(1+\eta)^\alpha}{\delta (C_\epsilon^*(\bm\mu))^\alpha}\right) \right].\]
This last upper bound yields, for every $\eta>0$ and $\epsilon>0$, 
\[\liminf_{\delta \rightarrow 0} \frac{\bE_{\bm \mu}[\tau_\delta]}{\log(1/\delta)} \leq \frac{\alpha (1+\eta)}{C_\epsilon^*(\bm \mu)}.\]
Letting $\eta$ and $\epsilon$ go to zero and using that, by continuity of $g$ and by definition of $w^*$,
\[\lim_{\epsilon \rightarrow 0} C_\epsilon^*(\bm \mu) = T^*(\bm \mu)^{-1}\]
yields 
\[\liminf_{\delta \rightarrow 0} \frac{\bE_{\bm \mu}[\tau_\delta]}{\log(1/\delta)} \leq \alpha T^*(\bm \mu)\;.\]

\subsection{Proof of Lemma~\ref{lem:concSimple}.} 
\[\bP\left(\cE_T^c\right) \leq \sum_{t=h(T)}^{T} \bP\left(\hat{\bm\mu}(t) \notin \cI_\epsilon \right) = \sum_{t=h(T)}^{T}\sum_{a=1}^{K} \left[\bP\left(\hat{\mu}_a(t) \leq \mu_a  - \xi\right) +  \bP\left(\hat{\mu}_a(t) \geq \mu_a + \xi\right)\right]\;. \]
Let $T$ be such that $h(T) \geq K^2$. Then for $t \geq h(T)$ one has $N_a(t) \geq (\sqrt{t} - K/2)_+ -1 \geq \sqrt{t} - K$ for every arm $a$.  Let $\hat{\mu}_{a,s}$ be the empirical mean of the first $s$ rewards from arm $a$ (such that $\hat{\mu}_a(t) = \hat{\mu}_{a,N_a(t)}$). Using a union bound and Chernoff inequality, one can write
\begin{eqnarray*}
 \bP\left(\hat{\mu}_a(t) \leq \mu_a - \xi\right) & = &  \bP\left(\hat{\mu}_a(t) \leq \mu_a - \xi, N_a(t) \geq \sqrt{t}\right)  \leq  \sum_{s=\sqrt{t} - K}^t \bP\left(\hat{\mu}_{a,s} \leq \mu_a - \xi\right) \\
 & \leq & \sum_{s=\sqrt{t}-K}^t \exp(-s d(\mu_a - \xi,\mu_a)) \leq \frac{e^{-(\sqrt{t}-K)d(\mu_a - \xi,\mu_a)}}{1 - e^{-d(\mu_a - \xi,\mu_a)}}\;.
\end{eqnarray*}
Similarly, one can prove that 
\begin{eqnarray*}
 \bP\left(\hat{\mu}_a(t) \geq \mu_a + \xi\right)\leq \frac{e^{-(\sqrt{t}-K)d(\mu_a + \xi,\mu_a)}}{1 - e^{-d(\mu_a + \xi,\mu_a)}}.
\end{eqnarray*}
Finally, letting 
\[C = \min_a\left(d(\mu_a - \xi,\mu_a)\wedge d(\mu_a+ \xi,\mu_a)\right) \  \ \text{and}  \ \ B = \sum_{a=1}^{K}\left( \frac{e^{Kd(\mu_a-\xi,\mu_a)}}{1 - e^{-d(\mu_a - \xi,\mu_a)}} + \frac{e^{Kd(\mu_a+\xi,\mu_a)}}{1 - e^{-d(\mu_a + \xi,\mu_a)}}\right)\;,\]
one obtains 
\[\bP\left(\cE_T^c\right) \leq \sum_{t=h(T)}^{T} B\exp(-\sqrt{t} C) \leq  B T \exp(-\sqrt{h(T)} C) =  B T \exp(-C T^{1/8})\;. \]

\end{document}